\documentclass[a4paper, 12pt]{article}
\usepackage{amssymb, amsthm,amsmath}
\usepackage{enumerate}
\usepackage{stmaryrd}
\usepackage{mathrsfs}
\usepackage{hyperref}
\usepackage[all]{xy}

\numberwithin{equation}{section}
\newtheorem{thm}[equation]{Theorem}
\newtheorem{lem}[equation]{Lemma}
\newtheorem{prop}[equation]{Proposition}
\newtheorem{cor}[equation]{Corollary}

\theoremstyle{definition}

\theoremstyle{remark}
\newtheorem{rem}[equation]{Remark}
\newtheorem{rems}[equation]{Remarks}

\newtheorem{example}[equation]{Example}

\newcommand{\im}{\operatorname{im}}
\newcommand{\coim}{\operatorname{coim}}
\newcommand{\Hom}{\operatorname{Hom}}
\newcommand{\id}{\operatorname{id}}
\newcommand{\coker}{\operatorname{coker}}
\newcommand{\Res}{\operatorname{Res}}
\newcommand{\Ind}{\operatorname{Ind}}
\newcommand{\Coind}{\operatorname{Coind}}
\newcommand{\Tot}{\operatorname{Tot}}
\newcommand{\Ext}{\operatorname{Ext}}
\newcommand{\cExt}{\operatorname{cExt}}
\newcommand{\rExt}{\operatorname{rExt}}
\newcommand{\Tor}{\operatorname{Tor}}
\newcommand{\cTor}{\operatorname{cTor}}
\newcommand{\rTor}{\operatorname{rTor}}
\newcommand{\FP}{\operatorname{FP}}

\begin{document}

\title{Continuous cohomology and homology\\ of profinite groups}
\author{Marco Boggi \and Ged Corob Cook}
\date{}
\maketitle

\begin{abstract}
We develop cohomological and homological theories for a profinite group $G$ with coefficients in the Pontryagin dual categories of pro-discrete and ind-profinite $G$-modules, respectively. The standard results of group (co)homology hold for this theory: we prove versions of the Universal Coefficient Theorem, the Lyndon-Hochschild-Serre spectral sequence and Shapiro's Lemma.
\end{abstract}

\section*{Introduction}

Cohomology groups $H^n(G,M)$ can be studied for profinite groups $G$ in much the same way as abstract groups. The coefficients $M$ will lie in some category of topological modules, but it is not clear what the right category is. The classical solution is to allow only discrete modules, in which case $H^n(G,M)$ is discrete: see \cite{R-Z} for this approach. For many applications, it is useful to take $M$ to be a profinite $G$-module. A cohomology theory allowing discrete and profinite coefficients is developed in \cite{S-W} when $G$ is of type $\FP_\infty$, but for arbitrary profinite groups there has not previously been a satisfactory definition of cohomology with profinite coefficients. A difficulty is that the category of profinite $G$-modules does not have enough injectives.

We define the cohomology of a profinite group with coefficients in the category of pro-discrete $\hat{\mathbb{Z}}\llbracket G\rrbracket$-modules, $PD(\hat{\mathbb{Z}}\llbracket G\rrbracket)$. This category contains the discrete $\hat{\mathbb{Z}}\llbracket G\rrbracket$-modules and the second-countable profinite $\hat{\mathbb{Z}}\llbracket G\rrbracket$-modules; when $G$ itself is second-countable, this is sufficient for many applications.

$PD(\hat{\mathbb{Z}}\llbracket G\rrbracket)$ is not an abelian category: instead it is quasi-abelian -- homological algebra over this generalisation is treated in detail in \cite{Prosmans} and \cite{Schneiders}, and we give an overview of the results we will need in Section \ref{qacs}. Working over the derived category, this allows us to define derived functors and study their properties: these functors exist because $PD(\hat{\mathbb{Z}}\llbracket G\rrbracket)$ has enough injectives. The resulting cohomology theory does not take values in a module category, but rather in the heart of a canonical t-structure on the derived category, $\mathcal{RH}(PD(\hat{\mathbb{Z}}))$, in which $PD(\hat{\mathbb{Z}})$ is a coreflective subcategory.

The most important result of this theory is that it allows us to prove a Lyndon-Hochschild-Serre spectral sequence for profinite groups with profinite coefficients. This has not been possible in previous formulations of profinite cohomology, and should allow the application of a wide range of techniques from abstract group cohomology to the study of profinite groups. A good example of this is the use of the spectral sequence to give a partial answer to a conjecture of Kropholler's, \cite[Open Question 6.12.1]{R-Z}, in a paper by the second author \cite{CC}.

We also define a homology theory for profinite groups which extends the category of coefficient modules to the ind-profinite $G$-modules. As in previous expositions, this is entirely dual to the cohomology theory.

Finally, in Section \ref{compare} we compare this theory to previous cohomology theories for profinite groups. It is naturally isomorphic to the classical cohomology of profinite groups with discrete coefficients, and to the Symonds-Weigel theory for profinite modules of type $\FP_\infty$ with profinite coefficients. We also define a continuous cochain cohomology, constructed by considering only the continuous $G$-maps from the standard bar resolution of a topological group $G$ to a topological $G$-module $M$, with the compact-open topology, and taking its cohomology; the comparison here is more nuanced, but we show that in certain circumstances these cohomology groups can be recovered from ours.

To clarify some terminology: it is common to refer to groups, modules, and so on without a topology as discrete. However, this creates an ambiguity in this situation. For a profinite ring $R$, there are $R$-modules $M$ without a topology such that giving $M$ the discrete topology creates a topological group on which the $R$-action $R \times M \to M$ is not continuous. Therefore a discrete module will mean one for which the $R$-action is continuous, and we will call algebraic objects without a topology \emph{abstract}.

\section{Ind-Profinite Modules}
\label{ip}

We say a topological space $X$ is \emph{ind-profinite} if there is an injective sequence of subspaces $X_i$, $i \in \mathbb{N}$, whose union is $X$, such that each $X_i$ is profinite and $X$ has the colimit topology with respect to the inclusions $X_i \rightarrow X$. That is, $X = \varinjlim_{IPSpace} X_i$. We write $IPSpace$ for the category of ind-profinite spaces and continuous maps.

\begin{prop}
\label{ipspace}
Given an ind-profinite space $X$ defined as the colimit of an injective sequence $\{X_i\}$ of profinite spaces, any compact subspace $K$ of $X$ is contained in some $X_i$.
\end{prop}
\begin{proof}
\cite[Proposition 1.1]{E-H} proves this under the additional assumption that the $X_i$ are profinite groups, but the proof does not use this.
\end{proof}

This shows that compact subspaces of $X$ are exactly the profinite subspaces, and that, if an ind-profinite space $X$ is defined as the colimit of a sequence $\{X_i\}$, then the $X_i$ are cofinal in the poset of compact subspaces of $X$. We call such a sequence a \emph{cofinal sequence for $X$}: any cofinal sequence of profinite subspaces defines $X$ up to homeomorphism.

A topological space $X$ is called \emph{compactly generated} if it satisfies the following condition: a subspace $U$ of $X$ is closed if and only if $U \cap K$ is closed in $K$ for every compact subspace $K$ of $X$. See \cite{Strickland} for background on such spaces. By the definition of the colimit topology, ind-profinite spaces are compactly generated. Indeed, a subspace $U$ of an ind-profinite space $X$ is closed if and only if $U \cap X_i$ is closed in $X_i$ for all $i$, if and only if $U \cap K$ is closed in $K$ for every compact subspace $K$ of $X$ by Proposition \ref{ipspace}.

\begin{lem}
\label{prodcoprod}
$IPSpace$ has finite products and coproducts.
\end{lem}
\begin{proof}
Given $X,Y \in IPSpace$ with cofinal sequences $\{X_i\},\{Y_i\}$, we can construct $X \sqcup Y$ using the cofinal sequence $\{X_i \sqcup Y_i\}$. However, it is not clear whether $X \times Y$ with the product topology is ind-profinite. Instead, thanks to Proposition \ref{ipspace}, the ind-profinite space $\varinjlim X_i \times Y_i$ is the product of $X$ and $Y$: it is easy to check that it satisfies the relevant universal property.
\end{proof}

Moreover, by the proposition, $\{X_i \times Y_i\}$ is cofinal in the poset of compact subspaces of $X \times Y$ (with the product topology), and hence $\varinjlim X_i \times Y_i$ is the \emph{$k$-ification} of $X \times Y$, or in other words it is the product of $X$ and $Y$ in the category of compactly generated spaces -- see \cite{Strickland} for details. So we will write $X \times_k Y$ for the product in $IPSpace$.

We say an abelian group $M$ equipped with an ind-profinite topology is an \emph{ind-profinite} abelian group if it satisfies the following condition: there is an injective sequence of profinite subgroups $M_i$, $i \in \mathbb{N}$, which is a cofinal sequence for the underlying space of $M$. It is easy to see that profinite groups and countable discrete torsion groups are ind-profinite. Moreover $\mathbb{Q}_p$ is ind-profinite via the cofinal sequence
\begin{equation*}
\mathbb{Z}_p \xrightarrow{\cdot p} \mathbb{Z}_p \xrightarrow{\cdot p} \cdots. \tag{$\ast$}
\end{equation*}

\begin{rem}
It is not obvious that ind-profinite abelian groups are topological groups. In fact, we see below that they are. But it is much easier to see that they are $k$-groups in the sense of \cite{Lamartin}: the multiplication map $M \times_k M = \varinjlim_{IPSpace} M_i \times M_i \rightarrow M$ is continuous by the definition of colimits. The $k$-group intuition will often be more useful.
\end{rem}

In the terminology of \cite{E-H} the ind-profinite abelian groups are just the abelian \emph{weakly profinite} groups. We recall some of the basic results of \cite{E-H}.

\begin{prop}
\label{ipgroup}
Suppose $M$ is an ind-profinite abelian group with cofinal sequence $\{M_i\}$.
\begin{enumerate}[(i)]
\item Any compact subspace of $M$ is contained in some $M_i$.
\item Closed subgroups $N$ of $M$ are ind-profinite, with cofinal sequence $N \cap M_i$.
\item Quotients of $M$ by closed subgroups $N$ are ind-profinite, with cofinal sequence $M_i/(N \cap M_i)$.
\item Ind-profinite abelian groups are topological groups.
\end{enumerate}
\end{prop}
\begin{proof}
\cite[Proposition 1.1, Proposition 1.2, Proposition 1.5]{E-H}
\end{proof}

As before, we call a sequence $\{M_i\}$ of profinite subgroups making $M$ into an ind-profinite group a cofinal sequence for $M$.

Suppose from now on that $R$ is a commutative profinite ring and $\Lambda$ is a profinite $R$-algebra.

\begin{rem}
We could define ind-profinite rings as colimits of injective sequences (indexed by $\mathbb{N}$) of profinite rings, and much of what follows does hold in some sense for such rings, but not much is lost by the restriction. In particular, it would be nice to use the machinery of ind-profinite rings to study $\mathbb{Q}_p$, but the sequence $(\ast)$ making $\mathbb{Q}_p$ into an ind-profinite abelian group does not make it into an ind-profinite ring because the maps are not maps of rings.
\end{rem}

We say that $M$ is a left $\Lambda$-$k$-module if $M$ is a $k$-group equipped with a continuous map $\Lambda \times_k M \rightarrow M$. A $\Lambda$-$k$-module homomorphism $M \rightarrow N$ is a continuous map which is a homomorphism of the underlying abstract $\Lambda$-modules. Because $\Lambda$ is profinite, $\Lambda \times_k M = \Lambda \times M$, so $\Lambda \times M \rightarrow M$ is continuous. Hence if $M$ is a topological group (that is, if multiplication $M \times M \rightarrow M$ is continuous) then it is a topological $\Lambda$-module.

We say that a left $\Lambda$-$k$-module $M$ equipped with an ind-profinite topology is a left \emph{ind-profinite $\Lambda$-module} if there is an injective sequence of profinite submodules $M_i$, $i \in \mathbb{N}$, which is a cofinal sequence for the underlying space of $M$. So countable discrete $\Lambda$-modules are ind-profinite, because finitely generated discrete $\Lambda$-modules are finite, and so are profinite $\Lambda$-modules. In particular $\Lambda$, with left-multiplication, is an ind-profinite $\Lambda$-module. Note that, since profinite $\hat{\mathbb{Z}}$-modules are the same as profinite abelian groups, ind-profinite $\hat{\mathbb{Z}}$-modules are the same as ind-profinite abelian groups.

Then we immediately get the following.

\begin{cor}
\phantomsection
\label{ipmod}
Suppose $M$ is an ind-profinite $\Lambda$-module with cofinal sequence $\{M_i\}$.
\begin{enumerate}[(i)]
\item Any compact subspace of $M$ is contained in some $M_i$.
\item Closed submodules $N$ of $M$ are ind-profinite, with cofinal sequence $N \cap M_i$.
\item Quotients of $M$ by closed submodules $N$ are ind-profinite, with cofinal sequence $M_i/(N \cap M_i)$.
\item Ind-profinite $\Lambda$-modules are topological $\Lambda$-modules.
\end{enumerate}
\end{cor}

As before, we call a sequence $\{M_i\}$ of profinite submodules making $M$ into an ind-profinite $\Lambda$-module a cofinal sequence for $M$.

\begin{lem}
\label{funsystem}
Ind-profinite $\Lambda$-modules have a fundamental system of neighbourhoods of $0$ consisting of open submodules. Hence such modules are Hausdorff and totally disconnected.
\end{lem}
\begin{proof}
Suppose $M$ has cofinal sequence $M_i$, and suppose $U \subseteq M$ is open, with $0 \in U$; by definition, $U \cap M_i$ is open in $M_i$ for all $i$. Profinite modules have a fundamental system of neighbourhoods of $0$ consisting of open submodules, by \cite[Lemma 5.1.1]{R-Z}, so we can pick an open submodule $N_0$ of $M_0$ such that $N_0 \subseteq U \cap M_0$. Now we proceed inductively: given an open submodule $N_i$ of $M_i$ such that $N_i \subseteq U \cap M_i$, let $f$ be the quotient map $M \rightarrow M/N_i$. Then $f(U)$ is open in $M/N_i$ by \cite[Proposition 1.3]{E-H}, so $f(U) \cap M_{i+1}/N_i$ is open in $M_{i+1}/N_i$. Pick an open submodule of $M_{i+1}/N_i$ which is contained in $f(U) \cap M_{i+1}/N_i$ and write $N_{i+1}$ for its preimage in $M_{i+1}$. Finally, let $N$ be the submodule of $M$ with cofinal sequence $\{N_i\}$: $N$ is open and $N \subseteq U$, as required.
\end{proof}

Write $IP(\Lambda)$ for the category whose objects are left ind-profinite $\Lambda$-modules, and whose morphisms $M \rightarrow N$ are $\Lambda$-$k$-module homomorphisms. We will identify the category of right ind-profinite $\Lambda$-modules with $IP(\Lambda^{op})$ in the usual way. Given $M \in IP(\Lambda)$ and a submodule $M'$, write $\overline{M'}$ for the closure of $M'$ in $M$. Given $M,N \in IP(\Lambda)$, write $\Hom_\Lambda^{IP}(M,N)$ for the abstract $R$-module of morphisms $M \rightarrow N$: this makes $\Hom_\Lambda^{IP}(-,-)$ into a functor $IP(\Lambda)^{op} \times IP(\Lambda) \rightarrow Mod(R)$ in the usual way, where $Mod(R)$ is the category of abstract $R$-modules and $R$-module homomorphisms.

\begin{prop}
$IP(\Lambda)$ is an additive category with kernels and cokernels.
\end{prop}
\begin{proof}
The category is clearly pre-additive; the biproduct $M \oplus N$ is the biproduct of the underlying abstract modules, with the topology of $M \times_k N$. The existence of kernels and cokernels follows from Corollary \ref{ipmod}; the cokernel of $f: M \rightarrow N$ is $N/\overline{f(M)}$.
\end{proof}

\begin{rem}
\label{nabelian}
The category $IP(\Lambda)$ is not abelian in general. Consider the countable direct sum $\oplus_{\aleph_0} \mathbb{Z}/2\mathbb{Z}$, with the discrete topology, and the countable direct product $\prod_{\aleph_0} \mathbb{Z}/2\mathbb{Z}$, with the profinite topology. Both are ind-profinite $\hat{\mathbb{Z}}$-modules. There is a canonical injective map $i: \oplus \mathbb{Z}/2\mathbb{Z} \rightarrow \prod \mathbb{Z}/2\mathbb{Z}$, but $i(\oplus \mathbb{Z}/2\mathbb{Z})$ is not closed in $\prod \mathbb{Z}/2\mathbb{Z}$. Moreover, $\oplus \mathbb{Z}/2\mathbb{Z}$ is not homeomorphic to $i(\oplus \mathbb{Z}/2\mathbb{Z})$, with the subspace topology, because $i(\oplus \mathbb{Z}/2\mathbb{Z})$ is not discrete, by the construction of the product topology.
\end{rem}

Given a morphism $f: M \rightarrow N$ in a category with kernels and cokernels, we write $\coim(f)$ for $\coker(\ker(f))$, and $\im(f)$ for $\ker(\coker(f))$. That is, $\coim(f) = f(M)$, with the quotient topology coming from $M$, and $\im(f) = \overline{f(M)}$, with the subspace topology coming from $N$. In an abelian category, $\coim(f) = \im(f)$, but the preceding remark shows that this fails in $IP(\Lambda)$.

We say a morphism $f: M \rightarrow N$ in $IP(\Lambda)$ is \emph{strict} if $\coim(f) = \im(f)$. In particular strict epimorphisms are surjections. Note that if $M$ is profinite all morphisms $f: M \rightarrow N$ must be strict, because compact subspaces of Hausdorff spaces are closed, so that $\coim(f) \rightarrow \im(f)$ is a continuous bijection of compact Hausdorff spaces and hence a topological isomorphism.

\begin{prop}
\label{ipsection}
Morphisms $f: M \rightarrow N$ in $IP(\Lambda)$ such that $f(M)$ is a closed subset of $N$ have continuous sections. So $f$ is strict in this case, and in particular continuous bijections are isomorphisms.
\end{prop}
\begin{proof}
\cite[Proposition 1.6]{E-H}
\end{proof}

\begin{cor}[Canonical decomposition of morphisms]
\label{ipdecompose}
Every morphism $f: M \rightarrow N$ in $IP(\Lambda)$ can be uniquely written as the composition of a strict epimorphism, a bimorphism and a strict monomorphism. Moreover the bimorphism is an isomorphism if and only if $f$ is strict.
\end{cor}
\begin{proof}
The decomposition is the usual one $M \rightarrow \coim(f) \xrightarrow{g} \im(f) \rightarrow N$, for categories with kernels and cokernels. Clearly $\coim(f) = f(M) \rightarrow N$ is injective, so $g$ is too, and hence $g$ is monic. Also the set-theoretic image of $M \rightarrow \im(f)$ is dense, so the set-theoretic image of $g$ is too, and hence $g$ is epic. Then everything follows from Proposition \ref{ipsection}.
\end{proof}

Because $IP(\Lambda)$ is not abelian, it is not obvious what the right notion of exactness is. We will say that a chain complex $$\cdots \rightarrow L \xrightarrow{f} M \xrightarrow{g} N \rightarrow \cdots$$ is \emph{strict exact at $M$} if $\coim(f) = \ker(g)$. We say a chain complex is strict exact if it is strict exact at each $M$.

Despite the failure of our category to be abelian, we can prove the following Snake Lemma, which will be useful later.

\begin{lem}
\label{snake}
Suppose we have a commutative diagram in $IP(\Lambda)$ of the form
\[
\xymatrix{& L \ar[r] \ar[d]^{f} & M \ar[r]^{p} \ar[d]^{g} & N \ar[r] \ar[d]^{h} & 0 \\
0 \ar[r] & L' \ar[r]^{i} & M' \ar[r] & N' &,}
\]
such that the rows are strict exact at $M,N,L',M'$ and $f,g,h$ are strict. Then we have a strict exact sequence $$\ker(f) \rightarrow \ker(g) \rightarrow \ker(h) \xrightarrow{\partial} \coker(f) \rightarrow \coker(g) \rightarrow \coker(h).$$
\end{lem}
\begin{proof}
Note that kernels in $IP(\Lambda)$ are preserved by forgetting the topology, and so are cokernels of strict morphisms by Proposition \ref{ipsection}. So by forgetting the topology and working with abstract $\Lambda$-modules we get the sequence described above from the standard Snake Lemma for abstract modules, which is exact as a sequence of abstract modules. This implies that, if all the maps in the sequence are continuous, then they have closed set-theoretic image, and hence the sequence is strict by Proposition \ref{ipsection}. To see that $\partial$ is continuous, we construct it as a composite of continuous maps. Since $\coim(p) = N$, by Proposition \ref{ipsection} again $p$ has a continuous section $s_1: N \rightarrow M$, and similarly $i$ has a continuous section $s_2: \im(i) \rightarrow L'$. Then, as usual, $\partial = s_2gs_1$. The continuity of the other maps is clear.
\end{proof}

\begin{prop}
The category $IP(\Lambda)$ has countable colimits.
\end{prop}
\begin{proof}
We show first that $IP(\Lambda)$ has countable direct sums. Given a countable collection $\{M_n: n \in \mathbb{N}\}$ of ind-profinite $\Lambda$-modules, write $\{M_{n,i}: i \in \mathbb{N}\}$, for each $n$, for a cofinal sequence for $M_n$. Now consider the injective sequence $\{N_n\}$ given by $N_n = \prod_{i=1}^n M_{i, n+1-i}$: each $N_n$ is a profinite $\Lambda$-module, so the sequence defines an ind-profinite $\Lambda$-module $N$. It is easy to check that the underlying abstract module of $N$ is $\bigoplus_n M_n$, that each canonical map $M_n \rightarrow N$ is continuous, and that any collection of continuous homomorphisms $M_n \rightarrow P$ in $IP(\Lambda)$ induces a continuous $N \rightarrow P$.

Now suppose we have a countable diagram $\{M_n\}$ in $IP(\Lambda)$. Write $S$ for the closed submodule of $\bigoplus M_n$ generated (topologically) by the elements with $j$th component $-x$, $k$th component $f(x)$ and all other components $0$, for all maps $f: M_j \rightarrow M_k$ in the diagram and all $x \in M_j$. By standard arguments, $(\bigoplus M_n)/S$, with the quotient topology, is the colimit of the diagram.
\end{proof}

\begin{rem}
\label{ipoplusexact}
We get from this construction that, given a countable collection of short strict exact sequences $$0 \rightarrow L_n \rightarrow M_n \rightarrow N_n \rightarrow 0$$ in $IP(\Lambda)$, their direct sum $$0 \rightarrow \bigoplus L_n \rightarrow \bigoplus M_n \rightarrow \bigoplus N_n \rightarrow 0$$ is strict exact by Proposition \ref{ipsection}, because the sequence of underlying modules is exact. So direct sums preserve kernels and cokernels, and in particular direct sums preserve strict maps, because given a countable collection of strict maps $\{f_n\}$ in $IP(\Lambda)$,
\begin{align*}
\coim(\bigoplus f_n) &= \coker(\ker(\bigoplus f_n)) = \bigoplus \coker(\ker(f_n)) \\
&= \bigoplus \ker(\coker(f_n)) = \ker(\coker(\bigoplus f_n)) = \im(\bigoplus f_n).
\end{align*}
\end{rem}

\begin{lem}
\phantomsection
\label{ipprojinj}
\begin{enumerate}[(i)]
\item For $M,N \in IP(\Lambda)$, let $\{M_i\},\{N_j\}$ cofinal sequences of $M$ and $N$, respectively, 
$\Hom_\Lambda^{IP}(M,N) = \varprojlim_i \varinjlim_j \Hom_\Lambda^{IP}(M_i,N_j)$, in the category of $R$-modules.
\item Given $X \in IPSpace$ with a cofinal sequence $\{X_i\}$ and $N \in IP(\Lambda)$ with cofinal sequence $\{N_j\}$, write $C(X,N)$ for the $R$-module of continuous maps $X \rightarrow N$. Then $C(X,N) = \varprojlim_i \varinjlim_j C(X_i,N_j)$.
\end{enumerate}
\end{lem}
\begin{proof}
\begin{enumerate}[(i)]
\item Since $M = \varinjlim_{IP(\Lambda)} M_i$, we have that 
$$\Hom_\Lambda^{IP}(M,N) = \varprojlim \Hom_\Lambda^{IP}(M_i,N).$$ 
Since the $N_j$ are cofinal for $N$, every continuous map $M_i \rightarrow N$ factors through some $N_j$, so $\Hom_\Lambda^{IP}(M_i,N) = \varinjlim \Hom_\Lambda^{IP}(M_i,N_j)$.
\item Similarly.
\end{enumerate}
\end{proof}

Given $X \in IPSpace$ as before, define a module $FX \in IP(\Lambda)$ in the following way: let $FX_i$ be the free profinite $\Lambda$-module on $X_i$. The maps $X_i \rightarrow X_{i+1}$ induce maps $FX_i \rightarrow FX_{i+1}$ of profinite $\Lambda$-modules, and hence we get an ind-profinite $\Lambda$-module with cofinal sequence $\{FX_i\}$. Write $FX$ for this module, which we will call the \emph{free ind-profinite $\Lambda$-module on $X$}.

\begin{prop}
\label{ipfree}
Suppose $X \in IPSpace$ and $N \in IP(\Lambda)$. Then we have $\Hom_\Lambda^{IP}(FX,N) = C(X,N)$, naturally in $X$ and $N$.
\end{prop}
\begin{proof}
First recall that, by the definition of free profinite modules, there holds $\Hom_\Lambda^{IP}(FX,N) = C(X,N)$ when $X$ and $N$ are profinite. Then by Lemma \ref{ipprojinj}, $$\Hom_\Lambda^{IP}(FX,N) = \varprojlim_i \varinjlim_j \Hom_\Lambda^{IP}(FX_i,N_j) = \varprojlim_i \varinjlim_j C(X_i,N_j) = C(X,N).$$ The isomorphism is natural because $\Hom_\Lambda^{IP}(F-,-)$ and $C(-,-)$ are both bifunctors.
\end{proof}

We call $P \in IP(\Lambda)$ \emph{projective} if $$0 \rightarrow \Hom_\Lambda^{IP}(P,L) \rightarrow \Hom_\Lambda^{IP}(P,M) \rightarrow \Hom_\Lambda^{IP}(P,N) \rightarrow 0$$ is an exact sequence in $Mod(R)$ whenever $$0 \rightarrow L \rightarrow M \rightarrow N \rightarrow 0$$ is strict exact. We will say $IP(\Lambda)$ has \emph{enough projectives} if for every $M \in IP(\Lambda)$ there is a projective $P$ and a strict epimorphism $P \rightarrow M$.

\begin{cor}
\label{ipenough}
$IP(\Lambda)$ has enough projectives.
\end{cor}
\begin{proof}
By Proposition \ref{ipfree} and Proposition \ref{ipsection}, $FX$ is projective for all $X \in IPSpace$. So given $M \in IP(\Lambda)$, $FM$ has the required property: the identity $M \rightarrow M$ induces a canonical `evaluation map' $\varepsilon: FM \rightarrow M$, which is strict epic because it is a surjection.
\end{proof}

\begin{lem}
Projective modules in $IP(\Lambda)$ are summands of free ones.
\end{lem}
\begin{proof}
Given a projective $P \in IP(\Lambda)$, pick a free module $F$ and a strict epimorphism $f: F \rightarrow P$. By definition, the map $\Hom_\Lambda^{IP}(P,F) \xrightarrow{f^\ast} \Hom_\Lambda^{IP}(P,P)$ induced by $f$ is a surjection, so there is some morphism $g: P \rightarrow F$ such that $f^\ast(g) = gf = \id_P$. Then we get that the map $\ker(f) \oplus P \rightarrow F$ is a continuous bijection, and hence an isomorphism by Proposition \ref{ipsection}.
\end{proof}

\begin{rems}
\phantomsection
\label{projsummand}
\begin{enumerate}[(i)]
\item We can also define the class of \emph{strictly free modules} to be free ind-profinite modules on ind-profinite spaces $X$ which have the form of a disjoint union of profinite spaces $X_i$. By the universal properties of coproducts and free modules we immediately get $FX = \bigoplus FX_i$. Moreover, for every ind-profinite space $Y$ there is some $X$ of this form with a surjection $X \rightarrow Y$: given a cofinal sequence $\{Y_i\}$ in $Y$, let $X = \bigsqcup Y_i$, and the identity maps $Y_i \rightarrow Y_i$ induce the required map $X \rightarrow Y$. Then the same argument as before shows that projective modules in $IP(\Lambda)$ are summands of strictly free ones.
\item Note that a profinite module in $IP(\Lambda)$ is projective in $IP(\Lambda)$ if and only if it is projective in the category of profinite $\Lambda$-modules. Indeed, Proposition \ref{ipfree} shows that free profinite modules are projective in $IP(\Lambda)$, and the rest follows.
\end{enumerate}
\end{rems}

\section{Pro-Discrete Modules}
\label{pd}

Write $PD(\Lambda)$ for the category of left \emph{pro-discrete} $\Lambda$-modules: the objects $M$ in this category are countable inverse limits, as topological $\Lambda$-modules, of discrete $\Lambda$-modules $M^i$, $i \in \mathbb{N}$; the morphisms are continuous $\Lambda$-module homomorphisms. So discrete torsion $\Lambda$-modules are pro-discrete, and so are second-countable profinite $\Lambda$-modules by \cite[Proposition 2.6.1, Lemma 5.1.1]{R-Z}, and in particular $\Lambda$, with left-multiplication, is a pro-discrete $\Lambda$-module if $\Lambda$ is second-countable. Moreover $\mathbb{Q}_p$ is a pro-discrete $\hat{\mathbb{Z}}$-module via the sequence $$\cdots \xrightarrow{\cdot p} \mathbb{Q}_p/\mathbb{Z}_p \xrightarrow{\cdot p} \mathbb{Q}_p/\mathbb{Z}_p.$$ We will identify the category of right pro-discrete $\Lambda$-modules with $PD(\Lambda^{op})$ in the usual way.

\begin{lem}
\label{1st}
Pro-discrete $\Lambda$-modules are first-countable.
\end{lem}
\begin{proof}
We can construct $M = \varprojlim M^i$ as a closed subspace of $\prod M^i$. Each $M^i$ is first-countable because it is discrete, and first-countability is closed under countable products and subspaces.
\end{proof}

\begin{rems}
\phantomsection
\label{pd=cg}
\begin{enumerate}[(i)]
\item This shows that $\Lambda$ itself can be regarded as a pro-discrete $\Lambda$-module if and only if it is first-countable, if and only if it is second-countable by \cite[Proposition 2.6.1]{R-Z}. Rings of interest are often second-countable; this class includes, for example, $\mathbb{Z}_p$, $\hat{\mathbb{Z}}$, $\mathbb{Q}_p$, and the completed group ring $R \llbracket G \rrbracket$ when $R$ and $G$ are second-countable.
\item Since first-countable spaces are always compactly generated by \cite[Proposition 1.6]{Strickland}, pro-discrete $\Lambda$-modules are compactly generated as topological spaces. In fact more is true. Given a pro-discrete $\Lambda$-module $M$ which is the inverse limit of a countable sequence $\{M^i\}$ of finite quotients, suppose $X$ is a compact subspace of $M$ and write $X^i$ for the image of $X$ in $M^i$. By compactness, each $X^i$ is finite. Let $N^i$ be the submodule of $M^i$ generated by $X^i$: because $X^i$ is finite, $\Lambda$ is compact and $M^i$ is discrete torsion, $N^i$ is finite. Hence $N = \varprojlim N^i$ is a profinite $\Lambda$-submodule of $M$ containing $X$. So pro-discrete modules $M$ are compactly generated by their profinite submodules $N$, in the sense that a subspace $U$ of $M$ is closed if and only if $U \cap N$ is closed in $N$ for all $N$.
\end{enumerate}
\end{rems}

\begin{lem}
\label{metrisable}
Pro-discrete $\Lambda$-modules are metrisable and complete.
\end{lem}
\begin{proof}
\cite[IX, Section 3.1, Proposition 1]{Bourbaki2} and the corollary to \cite[II, Section 3.5, Proposition 10]{Bourbaki}.
\end{proof}

In general, pro-discrete $\Lambda$-modules need not be second-countable, because for example $PD(\hat{\mathbb{Z}})$ contains uncountable discrete abelian groups. However, we have the following result.

\begin{lem}
Suppose a $\Lambda$-module $M$ has a topology which makes it pro-discrete and ind-profinite (as a $\Lambda$-module). Then $M$ is second-countable and locally compact.
\end{lem}
\begin{proof}
As an ind-profinite $\Lambda$-module, take a cofinal sequence of profinite submodules $M_i$. For any discrete quotient $N$ of $M$, the image of each $M_i$ in $N$ is compact and hence finite, and $N$ is the union of these images, so $N$ is countable. Then if $M$ is the inverse limit of a countable sequence of discrete quotients $M^j$, each $M^j$ is countable and $M$ can be identified with a closed subspace of $\prod M^j$, so $M$ is second-countable because second-countability is closed under countable products and subspaces. By Proposition \ref{metrisable}, $M$ is a Baire space, and hence by the Baire category theorem one of the $M_i$ must be open. The result follows.
\end{proof}

\begin{prop}
\phantomsection
\label{pdmod}
Suppose $M$ is a pro-discrete $\Lambda$-module which is the inverse limit of a sequence of discrete quotient modules $\{M^i\}$. Let $U^i = \ker(M \rightarrow M^i)$.
\begin{enumerate}[(i)]
\item The sequence $\{M^i\}$ is cofinal in the poset of all discrete quotient modules of $M$.
\item A closed submodule $N$ of $M$ is pro-discrete, with a cofinal sequence $\{N/(N \cap U^i)\}$.
\item Quotients of $M$ by closed submodules $N$ are pro-discrete, with cofinal sequence $\{M/(U^i+N)\}$.
\end{enumerate}
\end{prop}
\begin{proof}
\begin{enumerate}[(i)]
\item The $U^i$ form a basis of open neighbourhoods of $0$ in $M$, by \cite[Exercise 1.1.15]{R-Z}. Therefore, for any discrete quotient $D$ of $M$, the kernel of the quotient map $f: M \rightarrow D$ contains some $U^i$, so $f$ factors through $U^i$.
\item $M$ is complete, and hence $N$ is complete by \cite[II, Section 3.4, Proposition 8]{Bourbaki}. It is easy to check that $\{N \cap U^i\}$ is a fundamental system of neighbourhoods of the identity, so $N = \varprojlim N/(N \cap U^i)$ by \cite[III, Section 7.3, Proposition 2]{Bourbaki}. Also, since $M$ is metrisable, by \cite[IX, Section 3.1, Proposition 4]{Bourbaki2} $M/N$ is complete too. After checking that $(U^i+N)/N$ is a fundamental system of neighbourhoods of the identity in $M/N$, we get $M/N = \varprojlim M/(U^i+N)$ by applying \cite[III, Section 7.3, Proposition 2]{Bourbaki} again.
\end{enumerate}
\end{proof}

As a result of (i), we call $\{M^i\}$ a cofinal sequence for $M$.

As in $IP(\Lambda)$, it is clear from Proposition \ref{pdmod} that $PD(\Lambda)$ is an additive category with kernels and cokernels. 

Given $M,N \in PD(\Lambda)$, write $\Hom_\Lambda^{PD}(M,N)$ for the $R$-module of morphisms $M \rightarrow N$: this makes $\Hom_\Lambda^{PD}(-,-)$ into a functor $$PD(\Lambda)^{op} \times PD(\Lambda) \rightarrow Mod(R)$$ in the usual way. Note that the ind-profinite $\hat{\mathbb{Z}}$-modules in Remark \ref{nabelian} are also pro-discrete $\hat{\mathbb{Z}}$-modules, so the remark also shows that $PD(\Lambda)$ is not abelian in general.

As before, we say a morphism $f: M \rightarrow N$ in $PD(\Lambda)$ is \emph{strict} if $\coim(f) = \im(f)$. In particular strict epimorphisms are surjections. We say that a chain complex $$\cdots \rightarrow L \xrightarrow{f} M \xrightarrow{g} N \rightarrow \cdots$$ is \emph{strict exact at $M$} if $\coim(f) = \ker(g)$. We say a chain complex is strict exact if it is strict exact at each $M$.

\begin{rem}
In general, it is not clear whether a map $f: M \rightarrow N$ of pro-discrete modules with $f(M)$ closed in $N$ must be strict, as is the case for ind-profinite modules. However, we do have the following result.
\end{rem}

\begin{prop}
Let $f: M \rightarrow N$ be a morphism in $PD(\Lambda)$. Suppose that $M$ (and hence $\coim(f)$) is second-countable, and that the set-theoretic image $f(M)$ is closed in $N$. Then the continuous bijection $\coim(f) \rightarrow \im(f)$ is an isomorphism; in other words, $f$ is strict.
\end{prop}
\begin{proof}
\cite[Chapter 6, Problem R]{Kelley}
\end{proof}

As for ind-profinite modules, we can factorise morphisms in a canonical way.

\begin{cor}[Canonical decomposition of morphisms]
\label{pddecompose}
Every morphism $f: M \rightarrow N$ in $IP(\Lambda)$ can be uniquely written as the composition of a strict epimorphism, a bimorphism and a strict monomorphism. Moreover the bimorphism is an isomorphism if and only if the morphism is strict.
\end{cor}

\begin{rem}
\label{pdexact}
Suppose we have a short strict exact sequence $$0 \rightarrow L \xrightarrow{f} M \xrightarrow{g} N \rightarrow 0$$ in $PD(\Lambda)$. Pick a cofinal sequence $\{M^i\}$ for $M$. Then, as in Proposition \ref{pdmod}(ii), $L = \coim(f) = \im(f) = \varprojlim \im(\im(f) \rightarrow M^i)$, and similarly for $N$, so we can write the sequence as a surjective inverse limit of short (strict) exact sequences of discrete $\Lambda$-modules.

Conversely, suppose we have a surjective sequence of short (strict) exact sequences $$0 \rightarrow L^i \rightarrow M^i \rightarrow N^i \rightarrow 0$$ of discrete $\Lambda$-modules. Taking limits we get a sequence
\begin{equation*}
0 \rightarrow L \xrightarrow{f} M \xrightarrow{g} N \rightarrow 0 \tag{$\ast$}
\end{equation*}
of pro-discrete $\Lambda$-modules. It is easy to check that $\im(f) = \ker(g) = L = \coim(f)$, and $\coim(g) = \coker(f) = N = \im(g)$, so $f$ and $g$ are strict, and hence ($\ast$) is a short strict exact sequence.
\end{rem}

\begin{lem}
\label{pdprojinj}
Given $M,N \in PD(\Lambda)$, pick cofinal sequences $\{M^i\},\{N^j\}$ respectively. Then $\Hom_\Lambda^{PD}(M,N) = \varprojlim_j \varinjlim_i \Hom_\Lambda^{PD}(M^i,N^j)$, in the category of $R$-modules.
\end{lem}
\begin{proof}
Since $N = \varprojlim_{PD(\Lambda)} N^j$, we have by definition that $\Hom_\Lambda^{PD}(M,N) = \varprojlim \Hom_\Lambda^{PD}(M,N^j)$. Since the $M^i$ are cofinal for $M$, every continuous map $M \rightarrow N^j$ factors through some $M^i$, so $\Hom_\Lambda^{IP}(M,N^j) = \varinjlim \Hom_\Lambda^{IP}(M^i,N^j)$.
\end{proof}

We call $I \in PD(\Lambda)$ \emph{injective} if $$0 \rightarrow \Hom_\Lambda^{PD}(N,I) \rightarrow \Hom_\Lambda^{PD}(M,I) \rightarrow \Hom_\Lambda^{PD}(L,I) \rightarrow 0$$ is an exact sequence of $R$-modules whenever $$0 \rightarrow L \rightarrow M \rightarrow N \rightarrow 0$$ is strict exact.

\begin{lem}
\label{pdinjective}
Suppose that $I$ is a discrete $\Lambda$-module which is injective in the category of discrete $\Lambda$-modules. Then $I$ is injective in $PD(\Lambda)$.
\end{lem}
\begin{proof}
We know $\Hom_\Lambda^{PD}(-,I)$ is exact on discrete $\Lambda$-modules. Remark \ref{pdexact} shows that we can write short strict exact sequences of pro-discrete $\Lambda$-modules as surjective inverse limits of short exact sequences of discrete modules in $PD(\Lambda)$, and then, by injectivity, applying $\Hom_\Lambda^{PD}(-,I)$ gives a direct system of short exact sequences of $R$-modules; the exactness of such direct limits is well-known.
\end{proof}

In particular we get that $\mathbb{Q}/\mathbb{Z}$, with the discrete topology, is injective in $PD(\hat{\mathbb{Z}})$ -- it is injective among discrete $\hat{\mathbb{Z}}$-modules (i.e. torsion abelian groups) by Baer's lemma, because it is divisible (see \cite[2.3.1]{Weibel}).

Given $M \in IP(\Lambda)$, with a cofinal sequence $\{M_i\}$, and $N \in PD(\Lambda)$, with a cofinal sequence $\{N^j\}$, we can consider the continuous group homomorphisms $f: M \rightarrow N$ which are compatible with the $\Lambda$-action, i.e. such that $\lambda f(m) = f(\lambda m)$, for all $\lambda \in \Lambda, m \in M$. Consider the category $T(\Lambda)$ of topological $\Lambda$-modules and continuous $\Lambda$-module homomorphisms. We can consider $IP(\Lambda)$ and $PD(\Lambda)$ as full subcategories of $T(\Lambda)$, and observe that $M = \varinjlim_{T(\Lambda)} M_i$ and $N = \varprojlim_{T(\Lambda)} N^j$. We write $\Hom_\Lambda^T(M,N)$ for the $R$-module of morphisms $M \rightarrow N$ in $T(\Lambda)$. For the following lemma, this will denote an abstract $R$-module, after which we will define a topology on $\Hom_\Lambda^T(M,N)$ making it into a topological $R$-module.

\begin{lem}
\label{tprojlim}
As abstract $R$-modules, $\Hom_\Lambda^T(M,N) = \varprojlim_{i,j} \Hom_\Lambda^T(M_i,N^j)$.
\end{lem}

We may give each $\Hom_\Lambda^T(M_i,N^j)$ the discrete topology, which is also the compact-open topology in this case. Then we make $\varprojlim \Hom_\Lambda^T(M_i,N^j)$ into a topological $R$-module by giving it the limit topology: giving $\Hom_\Lambda^T(M,N)$ this topology therefore makes it into a pro-discrete $R$-module. From now on, $\Hom_\Lambda^T(M,N)$ will be understood to have this topology. The topology thus constructed is well-defined because the $M_i$ are cofinal for $M$ and the $N^j$ cofinal for $N$. Moreover, given a morphism $M \rightarrow M'$ in $IP(\Lambda)$, this construction makes the induced map $\Hom_\Lambda^T(M',N) \rightarrow \Hom_\Lambda^T(M,N)$ continuous, and similarly in the second variable, so that $\Hom_\Lambda^T(-,-)$ becomes a functor $IP(\Lambda)^{op} \times PD(\Lambda) \rightarrow PD(R)$. Of course the case when $M$ and $N$ are right $\Lambda$-modules behaves in the same way; we may express this by treating $M,N$ as left $\Lambda^{op}$-modules and writing $\Hom_{\Lambda^{op}}^T(M,N)$ in this case.

More generally, given a chain complex $$\cdots \xrightarrow{d_1} M_1 \xrightarrow{d_0} M_0 \xrightarrow{d_{-1}} \cdots$$ in $IP(\Lambda)$ and a cochain complex $$\cdots \xrightarrow{d^{-1}} N^0 \xrightarrow{d^0} N^1 \xrightarrow{d^1} \cdots$$ in $PD(\Lambda)$, both bounded below, let us define the double cochain complex $\{\Hom_\Lambda^T(M_p,N^q)\}$ with the obvious horizontal maps, and with the vertical maps defined in the obvious way except that they are multiplied by $-1$ whenever $p$ is odd: this makes $\Tot(\Hom_\Lambda^T(M_p,N^q))$ into a cochain complex which we denote by $\Hom_\Lambda^T(M,N)$. Each term in the total complex is the sum of finitely many pro-discrete $R$-modules, because $M$ and $N$ are bounded below, so $\Hom_\Lambda^T(M,N)$ is a complex in $PD(R)$.

Suppose $\Theta,\Phi$ are profinite $R$-algebras. Then let $PD(\Theta-\Phi)$ be the category of pro-discrete $\Theta-\Phi$-bimodules and continuous $\Theta-\Phi$-homomorphisms. If $M$ is an ind-profinite $\Lambda-\Theta$-bimodule and $N$ is a pro-discrete $\Lambda-\Phi$-bimodule, one can make $\Hom_\Lambda^T(M,N)$ into a pro-discrete $\Theta-\Phi$-bimodule in the same way as in the abstract case. We leave the details to the reader.

\section{Pontryagin Duality}

\begin{lem}
\label{QZexact}
Suppose that $I$ is a discrete $\Lambda$-module which is injective in $PD(\Lambda)$. Then $\Hom_\Lambda^T(-,I)$ sends short strict exact sequences of ind-profinite $\Lambda$-modules to short strict exact sequences of pro-discrete $R$-modules.
\end{lem}
\begin{proof}
Proposition \ref{ipsection} shows that we can write short strict exact sequences of ind-profinite $\Lambda$-modules as injective direct limits of short exact sequences of profinite modules in $IP(\Lambda)$, and then \cite[Exercise 5.4.7(b)]{R-Z} shows that applying $\Hom_\Lambda^T(-,I)$ gives a surjective inverse system of short exact sequences of discrete $R$-modules; the inverse limit of these is strict exact by Remark \ref{pdexact}.
\end{proof}

In particular this applies when $I=\mathbb{Q}/\mathbb{Z}$, with the discrete topology, as a $\hat{\mathbb{Z}}$-module.

Consider $\mathbb{Q}/\mathbb{Z}$, with the discrete topology, as an ind-profinite abelian group. Given $M \in IP(\Lambda)$, with a cofinal sequence $\{M_i\}$, we can think of $M$ as an ind-profinite abelian group by forgetting the $\Lambda$-action; then $\{M_i\}$ becomes a cofinal sequence of profinite abelian groups for $M$. Now apply $\Hom_{\hat{\mathbb{Z}}}^T(-,\mathbb{Q}/\mathbb{Z})$ to get a pro-discrete abelian group. We can endow each $\Hom_{\hat{\mathbb{Z}}}^T(M_i,\mathbb{Q}/\mathbb{Z})$ with the structure of a right $\Lambda$-module, such that the $\Lambda$-action is continuous, by \cite[p.165]{R-Z}. Taking inverse limits, we can therefore make $\Hom_{\hat{\mathbb{Z}}}^T(M,\mathbb{Q}/\mathbb{Z})$ into a pro-discrete right $\Lambda$-module, which we denote by $M^\ast$. As before, ${}^\ast$ gives a contravariant functor $IP(\Lambda) \rightarrow PD(\Lambda^{op})$. Lemma \ref{QZexact} now has the following immediate consequence.

\begin{cor}
\label{aststrict}
The functor ${}^\ast: IP(\Lambda) \rightarrow PD(\Lambda^{op})$ maps short strict exact sequences to short strict exact sequences.
\end{cor}

Suppose instead that $M \in PD(\Lambda)$, with a cofinal sequence $\{M^i\}$. As before, we can think of $M$ as a pro-discrete abelian group by forgetting the $\Lambda$-action, and then $\{M^i\}$ is a cofinal sequence of discrete abelian groups. Recall that, as (abstract) $\hat{\mathbb{Z}}$-modules, $\Hom_{\hat{\mathbb{Z}}}^{PD}(M,\mathbb{Q}/\mathbb{Z}) \cong \varinjlim_i \Hom_{\hat{\mathbb{Z}}}^{PD}(M^i,\mathbb{Q}/\mathbb{Z})$. We can endow each $\Hom_{\hat{\mathbb{Z}}}^{PD}(M^i,\mathbb{Q}/\mathbb{Z})$ with the structure of a profinite right $\Lambda$-module, by \cite[p.165]{R-Z}. Taking direct limits, we then make $\Hom_{\hat{\mathbb{Z}}}^{PD}(M,\mathbb{Q}/\mathbb{Z})$ into an ind-profinite right $\Lambda$-module, which we denote by $M_\ast$, and in the same way as before ${}_\ast$ gives a functor $PD(\Lambda) \rightarrow IP(\Lambda^{op})$.

Note that ${}_\ast$ also maps short strict exact sequences to short strict exact sequences, by Lemma \ref{pdinjective} and Proposition \ref{ipsection}. Note too that both ${}^\ast$ and ${}_\ast$ send profinite modules to discrete modules and vice versa; on such modules they give the same result as the usual Pontryagin duality functor of \cite[Section 2.9]{R-Z}.

\begin{thm}[Pontryagin duality]
The composite functors $IP(\Lambda) \xrightarrow{-^\ast} PD(\Lambda^{op}) \xrightarrow{-_\ast} IP(\Lambda)$ and $PD(\Lambda) \xrightarrow{-_\ast} IP(\Lambda^{op}) \xrightarrow{-^\ast} PD(\Lambda)$ are naturally isomorphic to the identity, so that $IP(\Lambda)$ and $PD(\Lambda)$ are dually equivalent.
\end{thm}
\begin{proof}
We give a proof for ${}^\ast\circ{}_\ast$; the proof for ${}_\ast\circ{}^\ast$ is similar. Given $M \in IP(\Lambda)$ with a cofinal sequence $M_i$, by construction $(M^\ast)_\ast$ has cofinal sequence $(M_i^\ast)_\ast$. By \cite[p.165]{R-Z}, the functors ${}^\ast$ and ${}_\ast$ give a dual equivalence between the categories of profinite and discrete $\Lambda$-modules, so we have natural isomorphisms $M_i \rightarrow (M_i^\ast)_\ast$ for each $i$, and the result follows.
\end{proof}

From now on, by abuse of notation, we will follow convention by writing ${}^\ast$ for both the functors ${}^\ast$ and ${}_\ast$.

\begin{cor}
Pontryagin duality preserves the canonical decomposition of morphisms. More precisely, given a morphism $f: M \rightarrow N$ in $IP(\Lambda)$, $\im(f)^\ast = \coim(f^\ast)$ and $\im(f^\ast) = \coim(f)^\ast$. In particular, $f^\ast$ is strict if and only if $f$ is. Similarly for morphisms in $PD(\Lambda)$.
\end{cor}
\begin{proof}
This follows from Pontryagin duality and the duality between the definitions of $\im$ and $\coim$. For the final observation, note that, by Corollary \ref{ipdecompose} and Corollary \ref{pddecompose},
\begin{align*}
f^\ast \text{ is strict} &\Leftrightarrow \im(f^\ast) = \coim(f^\ast) \\
&\Leftrightarrow \im(f) = \coim(f) \\
&\Leftrightarrow f \text{ is strict.}
\end{align*}
\end{proof}

\begin{cor}
\phantomsection
\label{duals}
\begin{enumerate}[(i)]
\item $PD(\Lambda)$ has countable limits.
\item Direct products in $PD(\Lambda)$ preserve kernels and cokernels, and hence strict maps.
\item $PD(\Lambda)$ has \emph{enough injectives}: for every $M \in PD(\Lambda)$ there is an injective $I$ and a strict monomorphism $M \rightarrow I$. A discrete $\Lambda$-module $I$ is injective in $PD(\Lambda)$ if and only if it is injective in the category of discrete $\Lambda$-modules.
\item Every injective in $PD(\Lambda)$ is a summand of a \emph{strictly cofree} one, i.e. one whose Pontryagin dual is strictly free.
\item Countable products of strict exact sequences in $PD(\Lambda)$ are strict exact.
\item Let $P$ be a profinite $\Lambda$-module which is projective in $IP(\Lambda)$. Then the functor $\Hom_\Lambda^T(P,-)$ sends strict exact sequences of pro-discrete $\Lambda$-modules to strict exact sequences of pro-discrete $R$-modules.
\end{enumerate}
\end{cor}

\begin{example}
\label{Qast}
It is easy to check that $\hat{\mathbb{Z}}^\ast = \mathbb{Q/Z}$ and $\mathbb{Z}_p^\ast = \mathbb{Q}_p/\mathbb{Z}_p$. Then $$\mathbb{Q}_p^\ast = (\varinjlim (\mathbb{Z}_p \xrightarrow{\cdot p} \mathbb{Z}_p \xrightarrow{\cdot p} \cdots))^\ast = \varprojlim (\cdots \xrightarrow{\cdot p} \mathbb{Q}_p/\mathbb{Z}_p \xrightarrow{\cdot p} \mathbb{Q}_p/\mathbb{Z}_p) = \mathbb{Q}_p.$$
\end{example}

The topology defined on $M^\ast = \Hom_{\hat{\mathbb{Z}}}^T(M,\mathbb{Q}/\mathbb{Z})$ when $M$ is an ind-profinite $\Lambda$-module coincides with the compact-open topology, because the (discrete) topology on each $\Hom_{\hat{\mathbb{Z}}}^T(M_i,\mathbb{Q}/\mathbb{Z})$ is the compact-open topology and every compact subspace of $M$ is contained in some $M_i$ by Proposition \ref{ipspace}. Similarly, for a pro-discrete $\Lambda$-module $N$, every compact subspace of $N$ is contained in some profinite submodule $L$ by Remark \ref{pd=cg}(ii), and so the compact-open topology on $\Hom_{\hat{\mathbb{Z}}}^{PD}(N,\mathbb{Q}/\mathbb{Z})$ coincides with the limit topology on $\varprojlim_{T(\Lambda)} \Hom_{\hat{\mathbb{Z}}}^{PD}(L,\mathbb{Q}/\mathbb{Z})$, where the limit is taken over all profinite submodules of $N$ and each $\Hom_{\hat{\mathbb{Z}}}^{PD}(L,\mathbb{Q}/\mathbb{Z})$ is given the (discrete) compact-open topology.

\begin{prop}
\label{coast}
The compact-open topology on $\Hom_{\hat{\mathbb{Z}}}^{PD}(N,\mathbb{Q}/\mathbb{Z})$ coincides with the topology defined on $N^\ast$.
\end{prop}
\begin{proof}
By the preceding remarks, $\Hom_{\hat{\mathbb{Z}}}^{PD}(N,\mathbb{Q}/\mathbb{Z})$ with the compact-open topology is just $\varprojlim_{\text{profinite } L \leq N} L^\ast$. So the canonical map $N^\ast \rightarrow \varprojlim L^\ast$ is a continuous bijection; we need to check it is open. By Lemma \ref{funsystem}, it suffices to check this for open submodules $K$ of $N^\ast$. Because $K$ is open, $N^\ast/K$ is discrete, so $(N^\ast/K)^\ast$ is a profinite submodule of $N$. Therefore there is a canonical continuous map $\varprojlim L^\ast \rightarrow (N^\ast/K)^{\ast\ast} = N^\ast/K$, whose kernel is open because $N^\ast/K$ is discrete. This kernel is $K$, and the result follows.
\end{proof}

\begin{cor}
The topology on ind-profinite $\Lambda$-modules is complete, Hausdorff and totally disconnected.
\end{cor}
\begin{proof}
By Lemma \ref{funsystem} we just need to show the topology is complete. Proposition \ref{coast} shows that ind-profinite $\Lambda$-modules are the inverse limit of their discrete quotients, and hence that the topology on such modules is complete, by the corollary to \cite[II, Section 3.5, Proposition 10]{Bourbaki}.
\end{proof}

Moreover, given ind-profinite $\Lambda$-modules $M,N$, the product $M \times_k N$ is the inverse limit of discrete modules $M' \times_k N'$, where $M'$ and $N'$ are discrete quotients of $M$ and $N$ respectively. But $M' \times_k N' = M' \times N'$, because both are discrete, so $M \times_k N = \varprojlim M' \times N' = M \times N$, the product in the category of topological modules.

\begin{prop}
\label{projstrictexact}
Suppose that $P \in IP(\Lambda)$ is projective. Then $\Hom_\Lambda^T(P,-)$ sends strict exact sequences in $PD(\Lambda)$ to strict exact sequences in $PD(R)$.
\end{prop}
\begin{proof}
For $P$ profinite this is Corollary \ref{duals}(vi). For $P$ strictly free, $P = \bigoplus P_i$, we get $\Hom_\Lambda^T(P,-) = \prod \Hom_\Lambda^T(P_i,-)$, which sends strict exact sequences to strict exact sequences because $\prod$ and $\Hom_\Lambda^T(P_i,-)$ do. Now the result follows from Remark \ref{projsummand}.
\end{proof}

\begin{lem}
\label{NastMast}
$\Hom_\Lambda^T(M,N) = \Hom_{\Lambda^{op}}^T(N^\ast,M^\ast)$ for all $M \in IP(\Lambda), N \in PD(\Lambda)$, naturally in both variables.
\end{lem}
\begin{proof}
Think of $\Hom_\Lambda^T(M,N)$ and $\Hom_{\Lambda^{op}}^T(N^\ast,M^\ast)$ as abstract $R$-modules. Then, the functor 
$\Hom_\Lambda^T(-,\mathbb{Q/Z})$ induces maps 
$$\Hom_\Lambda^T(M,N)\! \xrightarrow{f_1}\! \Hom_\Lambda^T(N^\ast,M^\ast)\! \xrightarrow{f_2}\! \Hom_\Lambda^T(N^{\ast\ast},M^{\ast\ast}) 
\!\xrightarrow{f_3}\! \Hom_\Lambda^T(N^{\ast\ast\ast},M^{\ast\ast\ast})$$ 
such that the compositions $f_2f_1$ and $f_3f_2$ are isomorphisms, so $f_2$ is an isomorphism. In particular, this holds when $M$ is profinite and $N$ is discrete, in which case the topology on $\Hom_\Lambda^T(M,N)$ is discrete; so, taking cofinal sequences $M_i$ for $M$ and $N^j$ for $N$, we get $\Hom_\Lambda^T(M_i,N^j) = \Hom_{\Lambda^{op}}^T(N^{j\ast},M_i^\ast)$ as topological modules for each $i,j$, and the topologies on $\Hom_\Lambda^T(M,N)$ and $\Hom_{\Lambda^{op}}^T(N^\ast,M^\ast)$ are given by the inverse limits of these. Naturality is clear.
\end{proof}

\begin{cor}
\label{injstrictexact}
Suppose that $I \in PD(\Lambda)$ is injective. Then $\Hom_\Lambda^T(-,I)$ sends strict exact sequences in $IP(\Lambda)$ to strict exact sequences in $PD(R)$.
\end{cor}

\begin{prop}[Baer's Lemma]
Suppose $I \in PD(\Lambda)$ is discrete. Then $I$ is injective in $PD(\Lambda)$ if and only if, for every closed left ideal $J$ of $\Lambda$, every map $J \rightarrow I$ extends to a map $\Lambda \rightarrow I$.
\end{prop}
\begin{proof}
Think of $\Lambda$ and $J$ as objects of $PD(\Lambda)$. The condition is clearly necessary. To see it is sufficient, suppose we are given a strict monomorphism $f: M \rightarrow N$ in $PD(\Lambda)$ and a map $g: M \rightarrow I$. Because $I$ is discrete, $\ker(g)$ is open in $M$. Because $f$ is strict, we can therefore pick an open submodule $U$ of $N$ such that $\ker(g) = M \cap U$. So the problem reduces to the discrete case: it is enough to show that $M/\ker(g) \rightarrow I$ extends to a map $N/U \rightarrow I$. In this case, the proof for abstract modules, \cite[Baer's Criterion 2.3.1]{Weibel}, goes through unchanged.
\end{proof}

Therefore a discrete $\hat{\mathbb{Z}}$-module which is injective in $PD(\hat{\mathbb{Z}})$ is divisible. On the other hand, the discrete $\hat{\mathbb{Z}}$-modules are just the torsion abelian groups with the discrete topology. So, by the version of Baer's Lemma for abstract modules (\cite[Baer's Criterion 2.3.1]{Weibel}), divisible discrete $\hat{\mathbb{Z}}$-modules are injective in the category of discrete $\hat{\mathbb{Z}}$-modules, and hence injective in $PD(\hat{\mathbb{Z}})$ too by Corollary \ref{duals}(iii). So duality gives:

\begin{cor}
\phantomsection
\label{proj=tf}
\begin{enumerate}[(i)]
\item A discrete $\hat{\mathbb{Z}}$-module is injective in $PD(\hat{\mathbb{Z}})$ if and only if it is divisible.
\item A profinite $\hat{\mathbb{Z}}$-module is projective in $IP(\hat{\mathbb{Z}})$ if and only if it is torsion-free.
\end{enumerate}
\end{cor}
\begin{proof}
Being divisible and being torsion-free are Pontryagin dual by \cite[Theorem 2.9.12]{R-Z}.
\end{proof}

\begin{rem}
\label{Qpnotinjective}
On the other hand, $\mathbb{Q}_p$ is not injective in $PD(\hat{\mathbb{Z}})$ (and hence not projective in $IP(\hat{\mathbb{Z}})$ either), despite being divisible (respectively, torsion-free). Indeed, consider the monomorphism $$f: \mathbb{Q}_p \rightarrow \prod_{\mathbb{N}}\mathbb{Q}_p/\mathbb{Z}_p, x \mapsto (x,x/p,x/p^2,\ldots),$$ which is strict because its dual $$f^\ast: \bigoplus_{\mathbb{N}}\mathbb{Z}_p \rightarrow \mathbb{Q}_p, (x_0, x_1, \ldots) \mapsto \sum_n x_n/p^n$$ is surjective and hence strict by Proposition \ref{ipsection}. Suppose $\mathbb{Q}_p$ is injective, so that $f$ splits; the map $g$ splitting it must send the torsion elements of $\prod_{\mathbb{N}}\mathbb{Q}_p/\mathbb{Z}_p$ to $0$ because $\mathbb{Q}_p$ is torsion-free. But the torsion elements contain $\bigoplus_{\mathbb{N}}\mathbb{Q}_p/\mathbb{Z}_p$, so they are dense in $\prod_{\mathbb{N}}\mathbb{Q}_p/\mathbb{Z}_p$ and hence $g=0$, giving a contradiction.
\end{rem}

Finally, we recall the definition of \emph{quasi-abelian} categories from \cite[Definition 1.1.3]{Schneiders}. Suppose that $\mathcal{E}$ is an additive category with kernels and cokernels. Now $f$ induces a unique canonical map $g: \coim(f) \rightarrow \im(f)$ such that $f$ factors as $$A \rightarrow \coim(f) \xrightarrow{g} \im(f) \rightarrow B,$$ and if $g$ is an isomorphism we say $f$ is \emph{strict}. We say $\mathcal{E}$ is a \emph{quasi-abelian} category if it satisfies the following two conditions:
\begin{enumerate}[(i)]
\item[(QA)] in any pull-back square
\[
\xymatrix{A' \ar[r]^{f'} \ar[d] & B' \ar[d] \\
A \ar[r]^{f} & B,}
\]
if $f$ is strict epic then so is $f'$;
\item[(QA$^\ast$)] in any push-out square
\[
\xymatrix{A \ar[r]^{f} \ar[d] & B \ar[d] \\
A' \ar[r]^{f'} & B',}
\]
if $f$ is strict monic then so is $f'$.
\end{enumerate}

$IP(\Lambda)$ satisfies axiom (QA) because forgetting the topology preserves pull-backs, and $Mod(\Lambda)$ satisfies (QA), so pull-backs of surjections are surjections. Recall by Remark \ref{pd=cg}(ii) that pro-discrete modules are compactly generated; hence $PD(\Lambda)$ satisfies (QA) by \cite[Proposition 2.36]{Strickland}, since the forgetful functor to topological spaces preserves pull-backs. Then both categories satisfy axiom (QA$^\ast$) by duality, and we have:

\begin{prop}
$IP(\Lambda)$ and $PD(\Lambda)$ are quasi-abelian categories.
\end{prop}

Moreover, note that the definition of a strict morphism in a quasi-abelian category agrees with our use of the term in $IP(\Lambda)$ and $PD(\Lambda)$.

\section{Tensor Products}
\label{tp}

As in the abstract case, we can define tensor products of ind-profinite modules. Suppose $L \in IP(\Lambda^{op}), M \in IP(\Lambda), N \in IP(R)$. We call a continuous map $b: L \times_k M \rightarrow N$ \emph{bilinear} if the following conditions hold for all $l,l_1,l_2 \in L, m,m_1,m_2 \in M, \lambda \in \Lambda$:
\begin{enumerate}[(i)]
\item $b(l_1+l_2,m) = b(l_1,m)+b(l_2,m)$;
\item $b(l,m_1+m_2) = b(l,m_1)+b(l,m_2)$;
\item $b(l\lambda,m) = b(l,\lambda m)$.
\end{enumerate}
Then $T \in IP(R)$, together with a bilinear map $\theta: L \times_k M \rightarrow T$, is the tensor product of $L$ and $M$ if, for every $N \in IP(R)$ and every bilinear map $b: L \times_k M \rightarrow N$, there is a unique morphism $f: T \rightarrow N$ in $IP(R)$ such that $b = f \theta$.

If such a $T$ exists, it is clearly unique up to isomorphism, and then we write $L \hat{\otimes}_\Lambda M$ for the tensor product. To show the existence of $L \hat{\otimes}_\Lambda M$, we construct it directly: $b$ defines a morphism $b': F(L \times_k M) \rightarrow N$ in $IP(R)$, where $F(L \times_k M)$ is the free ind-profinite $R$-module on $L \times_k M$. From the bilinearity of $b$, we get that the $R$-submodule $K$ of $F(L \times_k M)$ generated by the elements $$(l_1+l_2,m)-(l_1,m)-(l_2,m), (l,m_1+m_2)-(l,m_1)-(l,m_2), (l\lambda,m)-(l,\lambda m)$$ for all $l,l_1,l_2 \in L, m,m_1,m_2 \in M, \lambda \in \Lambda$ is mapped to $0$ by $b'$. From the continuity of $b'$ we get that its closure $\bar{K}$ is mapped to $0$ too. Thus $b'$ induces a morphism $b'': F(L \times_k M)/\bar{K} \rightarrow N$. Then it is not hard to check that $F(L \times_k M)/\bar{K}$, together with $b''$, satisfies the universal property of the tensor product.

\begin{prop}
\phantomsection
\label{tensorprops}
\begin{enumerate}[(i)]
\item $- \hat{\otimes}_\Lambda -$ is an additive bifunctor $IP(\Lambda^{op}) \times IP(\Lambda) \rightarrow IP(R)$.
\item There is an isomorphism $\Lambda \hat{\otimes}_\Lambda M = M$ for all $M \in IP(\Lambda)$, natural in $M$, and similarly $L \hat{\otimes}_\Lambda \Lambda = L$ naturally.
\item $L \hat{\otimes}_\Lambda M = M \hat{\otimes}_{\Lambda^{op}} L$, naturally in $L$ and $M$.
\item Given $L$ in $IP(\Lambda^{op})$ and $M$ in $IP(\Lambda)$, with cofinal sequences $\{L_i\}$ and $\{M_j\}$, there is an isomorphism $$L \hat{\otimes}_\Lambda M \cong \varinjlim_{IP(R)} (L_i \hat{\otimes}_\Lambda M_j).$$
\end{enumerate}
\end{prop}
\begin{proof}
(i) and (ii) follow from the universal property.
\begin{enumerate}[(i)]
\setcounter{enumi}{2}
\item Writing $\ast$ for the $\Lambda^{op}$-actions, a bilinear map $b_\Lambda: L \times M \rightarrow N$ (satisfying $b_\Lambda(l\lambda,m) = b_\Lambda(l,\lambda m)$) is the same thing as a bilinear map $b_{\Lambda^{op}}: M \times L \rightarrow N$ (satisfying $b_{\Lambda^{op}}(m,\lambda \ast l) = b_{\Lambda^{op}}(m \ast \lambda,l)$).
\item We have $L \times_k M = \varinjlim L_i \times M_j$ by Lemma \ref{prodcoprod}. By the universal property of the tensor product, the bilinear map $\varinjlim L_i \times M_j \rightarrow L \times_k M \rightarrow L \hat{\otimes}_\Lambda M$ factors through $f: \varinjlim L_i \hat{\otimes}_\Lambda M_j \rightarrow L \hat{\otimes}_\Lambda M$, and similarly the bilinear map $L \times_k M \rightarrow \varinjlim L_i \times M_j \rightarrow \varinjlim L_i \hat{\otimes}_\Lambda M_j$ factors through $g: L \hat{\otimes}_\Lambda M \rightarrow \varinjlim L_i \hat{\otimes}_\Lambda M_j$. By uniqueness, the compositions $fg$ and $gf$ are both identity maps, so the two sides are isomorphic.
\end{enumerate}
\end{proof}

More generally, given chain complexes $$\cdots \xrightarrow{d_1} L_1 \xrightarrow{d_0} L_0 \xrightarrow{d_{-1}} \cdots$$ in $IP(\Lambda^{op})$ and $$\cdots \xrightarrow{d'_1} M_1 \xrightarrow{d'_0} M_0 \xrightarrow{d'_{-1}} \cdots$$ in $IP(\Lambda)$, both bounded below, define the double chain complex $\{L_p \hat{\otimes}_\Lambda M_q\}$ with the obvious vertical maps, and with the horizontal maps defined in the obvious way except that they are multiplied by $-1$ whenever $q$ is odd: this makes $\Tot(L \hat{\otimes}_\Lambda M)$ into a chain complex which we denote by $L \hat{\otimes}_\Lambda M$. Each term in the total complex is the sum of finitely many ind-profinite $R$-modules, because $M$ and $N$ are bounded below, so $L \hat{\otimes}_\Lambda M$ is a complex in $IP(R)$.

Suppose from now on that $\Theta,\Phi,\Psi$ are profinite $R$-algebras. Then let $IP(\Theta-\Phi)$ be the category of ind-profinite $\Theta-\Phi$-bimodules and $\Theta-\Phi$-$k$-bimodule homomorphisms. We leave the details to the reader, after noting that an ind-profinite $R$-module $N$, with a left $\Theta$-action and a right $\Phi$-action which are continuous on profinite submodules, is an ind-profinite $\Theta-\Phi$-bimodule since we can replace a cofinal sequence $\{N_i\}$ of profinite $R$-modules with a cofinal sequence $\{\Theta \cdot N_i \cdot \Phi\}$ of profinite $\Theta-\Phi$-bimodules. If $L$ is an ind-profinite $\Theta-\Lambda$-bimodule and $M$ is an ind-profinite $\Lambda-\Phi$-bimodule, one can make $L \hat{\otimes}_\Lambda M$ into an ind-profinite $\Theta-\Phi$-bimodule in the same way as in the abstract case.

\begin{thm}[Adjunction isomorphism]
\label{adjoint}
Suppose $L \in IP(\Theta-\Lambda), M \in IP(\Lambda-\Phi), N \in PD(\Theta-\Psi)$. Then there is an isomorphism $$\Hom_\Theta^T(L \hat{\otimes}_\Lambda M,N) \cong \Hom_\Lambda^T(M,\Hom_\Theta^T(L,N))$$ in $PD(\Phi-\Psi)$, natural in $L,M,N$.
\end{thm}
\begin{proof}
Given cofinal sequences $\{L_i\},\{M_j\},\{N^k\}$ in $L,M,N$ respectively, we have natural isomorphisms 
$$\Hom_\Theta^T(L_i \hat{\otimes}_\Lambda M_j,N_k) \cong \Hom_\Lambda^T(M_j,\Hom_\Theta^T(L_i,N_k))$$ 
of discrete $\Phi-\Psi$-bimodules for each $i,j,k$ by \cite[Proposition 5.5.4(c)]{R-Z}. Then by Lemma \ref{tprojlim} we have
\begin{align*}
\Hom_\Theta^T(L \hat{\otimes}_\Lambda M,N) &\cong \varprojlim_{PD(\Phi-\Psi)} \Hom_\Theta^T(L_i \hat{\otimes}_\Lambda M_j,N_k) \\
&\cong \varprojlim_{PD(\Phi-\Psi)} \Hom_\Lambda^T(M_j,\Hom_\Theta^T(L_i,N_k)) \\
&\cong \Hom_\Lambda^T(M,\Hom_\Theta^T(L,N)).
\end{align*}
\end{proof}

It follows that $\Hom_\Lambda^T$ (considered as a co-/covariant bifunctor $IP(\Lambda)^{op} \times PD(\Lambda) \rightarrow PD(R)$) commutes with limits in both variables, and that $\hat{\otimes}_\Lambda$ commutes with colimits in both variables, by \cite[Theorem 2.6.10]{Weibel}.

If $L \in IP(\Theta-\Phi)$, Pontryagin duality gives $L^\ast$ the structure of a pro-discrete $\Phi-\Theta$-bimodule, and similarly with ind-profinite and pro-discrete switched.

\begin{cor}
\label{adjointast}
There is a natural isomorphism $$(L \hat{\otimes}_\Lambda M)^\ast \cong \Hom_\Lambda^T(M,L^\ast)$$ in $PD(\Phi-\Theta)$ for $L \in IP(\Theta-\Lambda), M \in IP(\Lambda-\Phi)$.
\end{cor}
\begin{proof}
Apply the theorem with $\Psi = \hat{\mathbb{Z}}$ and $N = \mathbb{Q}/\mathbb{Z}$.
\end{proof}

Properties proved about $\Hom_\Lambda$ in the past two sections carry over immediately to properties of $\hat{\otimes}_\Lambda$, using this natural isomorphism. Details are left to the reader.

Given a chain complex $M$ in $IP(\Lambda)$ and a cochain complex $N$ in $PD(\Lambda)$, both bounded below, if we apply ${}^\ast$ to the double complex with $(p,q)$th term $\Hom_\Lambda^T(M_p,N^q)$, we get a double complex with $(q,p)$th term $N^{q\ast} \hat{\otimes}_\Lambda M_p$ -- note that the indices are switched. This changes the sign convention used in forming $\Hom_\Lambda^T(M,N)$ into the one used in forming $N^\ast \hat{\otimes}_\Lambda M$, and so we have $\Hom_\Lambda^T(M,N)^\ast = N^\ast \hat{\otimes}_\Lambda M$ (because ${}^\ast$ commutes with finite direct sums).

\section{Derived Functors in Quasi-Abelian Categories}
\label{qacs}

We give a brief sketch of the machinery needed to derive functors in quasi-abelian categories. See \cite{Prosmans} and \cite{Schneiders} for details.

First a notational convention: in a chain complex $(A,d)$ in a quasi-abelian category, unless otherwise stated, $d_n$ will be the map $A_{n+1} \rightarrow A_n$. Dually, if $(A,d)$ is a cochain complex, $d^n$ will be the map $A^n \rightarrow A^{n+1}$.

Given a quasi-abelian category $\mathcal{E}$, let $\mathcal{K}(\mathcal{E})$ be the category whose objects are cochain complexes in $\mathcal{E}$ and whose morphisms are maps of cochain complexes up to homotopy; this makes $\mathcal{K}(\mathcal{E})$ into a triangulated category. Given a cochain complex $A$ in $\mathcal{E}$, we say $A$ is \emph{strict exact in degree $n$} if the map $d^{n-1}: A^{n-1} \rightarrow A^n$ is strict and $\im(d^{n-1}) = \ker(d^n)$. We say $A$ is \emph{strict exact} if it is strict exact in degree $n$ for all $n$. Then, writing $N(\mathcal{E})$ for the full subcategory of $\mathcal{K}(\mathcal{E})$ whose objects are strict exact, we get that $N(\mathcal{E})$ is a null system, so we can localise $\mathcal{K}(\mathcal{E})$ at $N(\mathcal{E})$ to get the derived category $\mathcal{D}(\mathcal{E})$. We also define $\mathcal{K}^+(\mathcal{E})$ to be the full subcategory of $\mathcal{K}(\mathcal{E})$ whose objects are bounded below, and $\mathcal{K}^-(\mathcal{E})$ to be the full subcategory whose objects are bounded above; we write $\mathcal{D}^+(\mathcal{E})$ and $\mathcal{D}^-(\mathcal{E})$ for their localisations, respectively. We say a map of complexes in $\mathcal{K}(\mathcal{E})$ is a \emph{strict quasi-isomorphism} if its cone is in $N(\mathcal{E})$.

Deriving functors in quasi-abelian categories uses the machinery of \emph{t-structures}. This can be thought of as giving a well-behaved cohomology functor to a triangulated category. For more detail on t-structures, see \cite[Section 1.3]{Prosmans}.

Given a triangulated category $\mathcal{T}$, with translation functor $T$, a t-structure on $\mathcal{T}$ is a pair $\mathcal{T}^{\leq 0}, \mathcal{T}^{\geq 0}$ of full subcategories of $\mathcal{T}$ satisfying the following conditions:
\begin{enumerate}[(i)]
\item $T(\mathcal{T}^{\leq 0}) \subseteq \mathcal{T}^{\leq 0}$ and $T^{-1}(\mathcal{T}^{\geq 0}) \subseteq \mathcal{T}^{\geq 0}$;
\item $\Hom_\mathcal{T}(X,Y) = 0$ for $X \in \mathcal{T}^{\leq 0}, Y \in T^{-1}(\mathcal{T}^{\geq 0})$;
\item for all $X \in \mathcal{T}$, there is a distinguished triangle $X_0 \rightarrow X \rightarrow X_1 \rightarrow$ with $X_0 \in \mathcal{T}^{\leq 0}, X_1 \in T^{-1}(\mathcal{T}^{\geq 0})$.
\end{enumerate}

It follows from this definition that, if $\mathcal{T}^{\leq 0}, \mathcal{T}^{\geq 0}$ is a t-structure on $\mathcal{T}$, there is a canonical functor $\tau^{\leq 0}: \mathcal{T} \rightarrow \mathcal{T}^{\leq 0}$ which is left adjoint to inclusion, and a canonical functor $\tau^{\geq 0}: \mathcal{T} \rightarrow \mathcal{T}^{\geq 0}$ which is right adjoint to inclusion. One can then define the \emph{heart} of the t-structure to be the full subcategory $\mathcal{T}^{\leq 0} \cap \mathcal{T}^{\geq 0}$, and the $0$th cohomology functor $$H^0: \mathcal{T} \rightarrow \mathcal{T}^{\leq 0} \cap \mathcal{T}^{\geq 0}$$ by $H^0 = \tau^{\geq 0}\tau^{\leq 0}$.

\begin{thm}
The heart of a t-structure on a triangulated category is an abelian category.
\end{thm}

There are two canonical t-structures on $\mathcal{D}(\mathcal{E})$, the left t-structure and the right t-structure, and correspondingly a left heart $\mathcal{LH(E)}$ and a right heart $\mathcal{RH(E)}$. The t-structures and hearts are dual to each other in the sense that there is a natural isomorphism between $\mathcal{LH(E)}$ and $\mathcal{RH}(\mathcal{E}^{op})$ (one can check that $\mathcal{E}^{op}$ is quasi-abelian), so we can restrict investigation to $\mathcal{LH(E)}$ without loss of generality.

Explicitly, the left t-structure on $\mathcal{D}(\mathcal{E})$ is given by taking $\mathcal{T}^{\leq 0}$ to be the complexes which are strict exact in all positive degrees, and $\mathcal{T}^{\geq 0}$ to be the complexes which are strict exact in all negative degrees. $\mathcal{LH(E)}$ is therefore the full subcategory of $\mathcal{D}(\mathcal{E})$ whose objects are strict exact in every degree except $0$; the $0$th cohomology functor $$LH^0: \mathcal{D}(\mathcal{E}) \rightarrow \mathcal{LH(E)}$$ is given by $$0 \rightarrow \coim(d^{-1}) \rightarrow \ker(d^0) \rightarrow 0.$$

Every object of $\mathcal{LH(E)}$ is isomorphic to a complex $$0 \rightarrow E^{-1} \xrightarrow{f} E^0 \rightarrow 0$$ of $\mathcal{E}$ with $E^0$ in degree $0$ and $f$ monic. Let $\mathcal{I}: \mathcal{E} \rightarrow \mathcal{LH(E)}$ be the functor given by $$E \mapsto (0 \rightarrow E \rightarrow 0)$$ with $E$ in degree $0$. Let $\mathcal{C}: \mathcal{LH(E)} \rightarrow \mathcal{E}$ be the functor given by $$(0 \rightarrow E^{-1} \xrightarrow{f} E^0 \rightarrow 0) \mapsto \coker(f).$$

\begin{prop}
\label{reflective}
$\mathcal{I}$ is fully faithful and right adjoint to $\mathcal{C}$. In particular, identifying $\mathcal{E}$ with its image under $\mathcal{I}$, we can think of $\mathcal{E}$ as a reflective subcategory of $\mathcal{LH(E)}$. Moreover, given a sequence $$0 \rightarrow L \rightarrow M \rightarrow N \rightarrow 0$$ in $\mathcal{E}$, its image under $\mathcal{I}$ is a short exact sequence in $\mathcal{LH(E)}$ if and only if the sequence is short strict exact in $\mathcal{E}$.
\end{prop}

The functor $\mathcal{I}$ induces a functor $\mathcal{D}(\mathcal{I}): \mathcal{D}(\mathcal{E}) \rightarrow \mathcal{D}(\mathcal{LH(E)})$.

\begin{prop}
\label{de=dlhe}
$\mathcal{D}(\mathcal{I})$ is an equivalence of categories which exchanges the left t-structure of $\mathcal{D(E)}$ with the standard t-structure of $\mathcal{D}(\mathcal{LH(E)})$. This induces equivalences $\mathcal{D(E)}^+ \rightarrow \mathcal{D}(\mathcal{LH(E)})^+$ and $\mathcal{D(E)}^- \rightarrow \mathcal{D}(\mathcal{LH(E)})^-$.
\end{prop}

Thus there are cohomological functors $LH^n: \mathcal{D}(\mathcal{E}) \rightarrow \mathcal{LH(E)}$, so that given any distinguished triangle in $\mathcal{D}(\mathcal{E})$ we get long exact sequences in $\mathcal{LH(E)}$. Given an object $(A,d) \in \mathcal{D}(\mathcal{E})$, $LH^n(A)$ is the complex $$0 \rightarrow \coim(d^{n-1}) \rightarrow \ker(d^n) \rightarrow 0$$ with $\ker(d^n)$ in degree $0$.

Everything for $\mathcal{RH(E)}$ is done dually, so in particular we get:

\begin{lem}
\label{rhoplh}
The functors $RH^n\colon \mathcal{D}(\mathcal{E}^{op}) \rightarrow \mathcal{RH}(\mathcal{E}^{op})$ are given by $$(LH^{-n})^{op}\colon \mathcal{D}(\mathcal{E})^{op} \rightarrow \mathcal{RH}(\mathcal{E})^{op}.$$
\end{lem}

As for $PD(\Lambda)$, we say an object $I$ of $\mathcal{E}$ is injective if, for any strict monomorphism $E \rightarrow E'$ in $\mathcal{E}$, any morphism $E \rightarrow I$ extends to a morphism $E' \rightarrow I$, and we say $\mathcal{E}$ has enough injectives if for every $E \in \mathcal{E}$ there is a strict monomorphism $E \rightarrow I$ for some injective $I$.

\begin{prop}
\label{injE=RHE}
The right heart $\mathcal{RH(E)}$ of $\mathcal{E}$ has enough injectives if and only if $\mathcal{E}$ does. An object $I \in \mathcal{E}$ is injective in $\mathcal{E}$ if and only if it is injective in $\mathcal{RH(E)}$.
\end{prop}

Suppose that $\mathcal{E}$ has enough injectives. Write $\mathcal{I}$ for the full subcategory of $\mathcal{E}$ whose objects are injective in $\mathcal{E}$.

\begin{prop}
Localisation at $N^+(\mathcal{I})$ gives an equivalence of categories $\mathcal{K}^+(\mathcal{I}) \rightarrow \mathcal{D}^+(\mathcal{E})$.
\end{prop}

We can now define derived functors in the same way as the abelian case. Suppose we are given an additive functor $F: \mathcal{E} \rightarrow \mathcal{E}'$ between quasi-abelian categories. Let $Q: \mathcal{K}^+(\mathcal{E}) \rightarrow \mathcal{D}^+(\mathcal{E})$ and $Q': \mathcal{K}^+(\mathcal{E}') \rightarrow \mathcal{D}^+(\mathcal{E}')$ be the canonical functors. Then the \emph{right derived functor} of $F$ is a triangulated functor $$RF: \mathcal{D}^+(\mathcal{E}) \rightarrow \mathcal{D}^+(\mathcal{E}')$$ (that is, a functor compatible with the triangulated structure) together with a natural transformation $$t: Q' \circ \mathcal{K}^+(F) \rightarrow RF \circ Q$$ satisfying the property that, given another triangulated functor $$G: \mathcal{D}^+(\mathcal{E}) \rightarrow \mathcal{D}^+(\mathcal{E}')$$ and a natural transformation $$g: Q' \circ \mathcal{K}^+(F) \rightarrow G \circ Q,$$ there is a unique natural transformation $h: RF \rightarrow G$ such that $g = (h \circ Q) t$. Clearly if $RF$ exists it is unique up to natural isomorphism.

Suppose we are given an additive functor $F: \mathcal{E} \rightarrow \mathcal{E}'$ between quasi-abelian categories, and suppose $\mathcal{E}$ has enough injectives. 

\begin{prop}
For $E \in \mathcal{K}^+(\mathcal{E})$ there is an $I \in \mathcal{K}^+(\mathcal{E})$ and a strict quasi-isomorphism $E \rightarrow I$ such that each $I^n$ is injective and each $E^n \rightarrow I^n$ is a strict monomorphism.
\end{prop}

We say such an $I$ is an \emph{injective resolution} of $E$.

\begin{prop}
\label{rf}
In the situation above, the right derived functor of $F$ exists and $RF(E) = \mathcal{K}^+(F)(I)$ for any injective resolution $I$ of $E$.
\end{prop}

We write $R^nF$ for the composition $RH^n \circ RF$.

\begin{rem}
\label{rrnotlr}
Since $RF$ is a triangulated functor, we could also define the cohomological functor $LH^n \circ RF$. The reason for using $RH^n \circ RF$ is Proposition \ref{injE=RHE}. Indeed, when $\mathcal{RH(E)}$ has enough injectives we may construct Cartan-Eilenberg resolutions in this category, and hence prove a Grothendieck spectral sequence, Theorem \ref{grothendieck} below. On the other hand it is not clear that such a spectral sequence holds for $LH^n \circ RF$, and in this sense $RH^n \circ RF$ is the `right' definition -- but see Lemma \ref{hom=ext}.
\end{rem}

The construction of derived functors generalises to the case of additive bifunctors $F: \mathcal{E} \times \mathcal{E}' \rightarrow \mathcal{E}''$ where $\mathcal{E}$ and $\mathcal{E}'$ have enough injectives: the right derived functor $$RF: \mathcal{D}^+(\mathcal{E}) \times \mathcal{D}^+(\mathcal{E}') \rightarrow \mathcal{D}^+(\mathcal{E}'')$$ exists and is given by $RF(E,E') = s\mathcal{K}^+(F)(I,I')$ where $I,I'$ are injective resolutions of $E,E'$ and $s\mathcal{K}^+(F)(I,I')$ is the total complex of the double complex $\{\mathcal{K}^+(F)(I^p,I'^q)\}_{pq}$ in which the vertical maps with $p$ odd are multiplied by $-1$.

Projectives are defined dually to injectives, left derived functors are defined dually to right derived ones, and if a quasi-abelian category $\mathcal{E}$ has enough projectives then an additive functor $F$ from $\mathcal{E}$ to another quasi-abelian category has a left derived functor $LF$ which can be calculated by taking projective resolutions, and we write $L_nF$ for $LH^{-n} \circ LF$. Similarly for bifunctors.

We state here, for future reference, some results on spectral sequences; see \cite[Chapter 5]{Weibel} for more details. All of the following results have dual versions obtained by passing to the opposite category, and we will use these dual results interchangeably with the originals. Suppose that $A = A^{pq}$ is a bounded below double cochain complex in $\mathcal{E}$, that is, there are only finitely many non-zero terms on each diagonal $n=p+q$, and the total complex $Tot(A)$ is bounded below. By Proposition \ref{de=dlhe}, we can equivalently think of $A$ as a bounded below double complex in the abelian category $\mathcal{RH(E)}$. Then we can use the usual spectral sequences for double complexes:

\begin{prop}
\label{doubless}
There are two bounded spectral sequences
\begin{equation*}
\begin{array}{rcr}
^IE_2^{pq} &=& RH_h^p RH_v^q (A) \\
^{II}E_2^{pq} &=& RH_v^p RH_h^q (A)
\end{array}
\Rightarrow RH^{p+q}\Tot(A),
\end{equation*}
naturally in $A$.
\end{prop}
\begin{proof}
\cite[Section 5.6]{Weibel}
\end{proof}

Suppose we are given an additive functor $F: \mathcal{E} \rightarrow \mathcal{E}'$ between quasi-abelian categories, and consider the case where $A \in \mathcal{D}^+(\mathcal{E})$. Suppose $\mathcal{E}$ has enough injectives, so that $\mathcal{RH(E)}$ does too. Thinking of $A$ as an object in $\mathcal{D}^+(\mathcal{RH(E)})$, we can take a bounded below Cartan-Eilenberg resolution $I$ of $A$. Then we can apply Proposition \ref{doubless} to the bounded below double complex $F(I)$ to get the following result.

\begin{prop}
\label{hyperhom}
There are two bounded spectral sequences
\begin{equation*}
\begin{array}{rcr}
^IE_2^{pq} &=& RH^p(R^qF(A)) \\
^{II}E_2^{pq} &=& (R^pF)(RH^q(A))
\end{array}
\Rightarrow R^{p+q}F(A),
\end{equation*}
naturally in $A$.
\end{prop}
\begin{proof}
\cite[Section 5.7]{Weibel}
\end{proof}

Suppose now that we are given additive functors $G: \mathcal{E} \rightarrow \mathcal{E}'$, $F: \mathcal{E}' \rightarrow \mathcal{E}''$ between quasi-abelian categories, where $\mathcal{E}$ and $\mathcal{E}'$ have enough injectives. Suppose $G$ sends injective objects of $\mathcal{E}$ to injective objects of $\mathcal{E}'$.

\begin{thm}[Grothendieck Spectral Sequence]
\label{grothendieck}
For $A \in \mathcal{D}^+(\mathcal{E})$ there is a natural isomorphism $R(FG)(A) \rightarrow (RF)(RG)(A)$ and a bounded spectral sequence $$^IE_2^{pq} = (R^pF)(R^qG(A)) \Rightarrow R^{p+q}(FG)(A),$$
naturally in $A$.
\end{thm}
\begin{proof}
Let $I$ be an injective resolution of $A$. There is a natural transformation $R(FG) \rightarrow (RF)(RG)$ by the universal property of derived functors; it is an isomorphism because, by hypothesis, each $G(I^n)$ is injective and hence $$(RF)(RG)(A) = F(G(I)) = R(FG)(A).$$ For the spectral sequence, apply Proposition \ref{hyperhom} with $A = G(I)$. We have $$^IE_2^{pq} = RH^p(R^qF(G(I))) \Rightarrow R^{p+q}F(G(I));$$ by the injectivity of the $G(I^n)$, $R^qF(G(I)) = 0$ for $q>0$, so the spectral sequence collapses to give $$R^{p+q}F(G(I)) \cong RH^p(FG(I)) = R^p(FG)(A).$$ On the other hand, $$^{II}E_2^{pq} = (R^pF)(RH^q(G(I))) = (R^pF)(R^qG(A))$$ and the result follows.
\end{proof}

We consider once more the case of an additive bifunctor $$F: \mathcal{E} \times \mathcal{E}' \rightarrow \mathcal{E}''$$ for $\mathcal{E}, \mathcal{E}'$ and $\mathcal{E}''$ quasi-abelian: this induces a triangulated functor $$\mathcal{K}^+(F): \mathcal{K}^+(\mathcal{E}) \times \mathcal{K}^+(\mathcal{E}') \rightarrow \mathcal{K}^+(\mathcal{E}''),$$ in the sense that a distinguished triangle in one of the variables, and a fixed object in the other, maps to a distinguished triangle in $\mathcal{K}^+(\mathcal{E}'')$. Hence for a fixed $A \in \mathcal{K^+(E)}$, $\mathcal{K}^+(F)$ restricts to a triangulated functor $\mathcal{K}^+(F)(A,-)$, and if $\mathcal{E}'$ has enough injectives we can derive this to get a triangulated functor $$R(F(A,-)): \mathcal{D}^+(\mathcal{E}') \rightarrow \mathcal{D}^+(\mathcal{E}'').$$ 

Maps $A \rightarrow A'$ in $\mathcal{K}^+(\mathcal{E})$ induce natural transformations $R(F(A,-)) \rightarrow R(F(A',-))$, so in fact we get a functor which we denote by $$R_2F: \mathcal{K}^+(\mathcal{E}) \times \mathcal{D}^+(\mathcal{E}') \rightarrow \mathcal{D}^+(\mathcal{E}'').$$ We know $R_2F$ is triangulated in the second variable, and it is triangulated in the first variable too because, given $B \in \mathcal{D}^+(\mathcal{E}')$ with an injective resolution $I$, $R_2F(-,B) = \mathcal{K}^+(F)(-,I)$ is a triangulated functor $\mathcal{K}^+(\mathcal{E}) \rightarrow \mathcal{D}^+(\mathcal{E}'')$. 

Similarly, we can define a triangulated functor $$R_1F: \mathcal{D}^+(\mathcal{E}) \times \mathcal{K}^+(\mathcal{E}') \rightarrow \mathcal{D}^+(\mathcal{E}'')$$ by deriving in the first variable, if $\mathcal{E}$ has enough injectives.

\begin{prop}
\phantomsection
\label{bifunctor}
\begin{enumerate}[(i)]
\item If $\mathcal{E}'$ has enough injectives and $F(-,J)\colon \mathcal{E} \rightarrow \mathcal{E}''$ is strict exact for $J$ injective, then $R_2F(-,B)$ sends quasi-isomorphisms to isomorphisms; that is, we can think of $R_2F$ as a functor $\mathcal{D}^+(\mathcal{E}) \times \mathcal{D}^+(\mathcal{E}') \rightarrow \mathcal{D}^+(\mathcal{E}'')$.
\item Suppose in addition that $\mathcal{E}$ has enough injectives. Then $R_2F$ is naturally isomorphic to $RF$.
\end{enumerate}
Similarly with the variables switched.
\end{prop}
\begin{proof}
\begin{enumerate}[(i)]
\item $R_2F(-,B) = \mathcal{K}^+(F)(-,I)$, for an injective resolution $I$ of $B$. Given a quasi-isomorphism $A \rightarrow A'$ in $\mathcal{K^+(E)}$, consider the map of double complexes $\mathcal{K}^+(F)(A,I) \rightarrow \mathcal{K}^+(F)(A',I)$ and apply Proposition \ref{doubless} to show that this map induces a quasi-isomorphism of the corresponding total complexes.
\item This holds by the same argument as (i), taking $A'$ to be an injective resolution of $A$.
\end{enumerate}
\end{proof}

\section{Derived Functors in \texorpdfstring{$IP(\Lambda)$}{IP(Lambda)} and \texorpdfstring{$PD(\Lambda)$}{PD(Lambda)}}

We now use the framework of Section \ref{qacs} to define derived functors in our categories of interest. Note first that the dual equivalence between $IP(\Lambda)$ and $PD(\Lambda)$ extends to dual equivalences between $\mathcal{D}^-(IP(\Lambda))$ and $\mathcal{D}^+(PD(\Lambda))$ given by applying the functor ${}^\ast$ to cochain complexes in these categories, by defining $(A^\ast)^n = (A^{-n})^\ast$ for a cochain complex $A$ in $PD(\Lambda)$, and similarly for the maps. We will also identify $\mathcal{D}^-(IP(\Lambda))$ with the category of chain complexes $A$ (localised over the strict quasi-isomorphisms) which are $0$ in negative degrees by setting $A_n = A^{-n}$. The Pontryagin duality extends to one between $\mathcal{LH}(IP(\Lambda))$ and $\mathcal{RH}(PD(\Lambda))$. Moreover, writing $RH^n$ and $LH^n$ for the $n$th cohomological functors $\mathcal{D}(PD(R)) \rightarrow \mathcal{RH}(PD(R))$ and $\mathcal{D}(IP(R^{op})) \rightarrow \mathcal{LH}(IP(R^{op}))$, respectively, the following is just a restatement of Lemma \ref{rhoplh}.

\begin{lem}
\label{heart*}
$LH^{-n}\circ{}^\ast = {}^\ast\circ RH^n$.
\end{lem}

Let $$R\Hom_\Lambda^T(-,-): \mathcal{D}^-(IP(\Lambda)) \times \mathcal{D}^+(PD(\Lambda)) \rightarrow \mathcal{D}^+(PD(R))$$ be the right derived functor of $\Hom_\Lambda^T(-,-): IP(\Lambda) \times PD(\Lambda) \rightarrow PD(R)$. By Proposition \ref{rf}, this exists because $IP(\Lambda)$ has enough projectives and $PD(\Lambda)$ has enough injectives, and $R\Hom_\Lambda^T(M,N)$ is given by $\Hom_\Lambda^T(P,I)$, where $P$ is a projective resolution of $M$ and $I$ is an injective resolution of $N$. Dually, let $$- \hat{\otimes}^L_\Lambda -: \mathcal{D}^-(IP(\Lambda^{op})) \times \mathcal{D}^-(IP(\Lambda)) \rightarrow \mathcal{D}^-(IP(R))$$ be the left derived functor of $- \hat{\otimes}_\Lambda -: IP(\Lambda^{op}) \times IP(\Lambda) \rightarrow IP(R)$. Then by Proposition \ref{rf} again $M \hat{\otimes}^L_\Lambda N$ is given by $P \hat{\otimes}_\Lambda Q$ where $P,Q$ are projective resolutions of $M,N$ respectively.

We also define $\Ext_\Lambda^n$ to be the composite
\begin{align*}
\mathcal{LH}(IP(\Lambda)) \times \mathcal{RH}(PD(\Lambda)) &\rightarrow \mathcal{D}^-(IP(\Lambda)) \times \mathcal{D}^+(PD(\Lambda)) \\
&\xrightarrow{R\Hom_\Lambda^T} \mathcal{D}^+(PD(R)) \\
&\xrightarrow{RH^n} \mathcal{RH}(PD(R))
\end{align*}
and $\Tor^\Lambda_n$ to be the composite
\begin{align*}
\mathcal{LH}(IP(\Lambda^{op})) \times \mathcal{LH}(IP(\Lambda)) &\rightarrow \mathcal{D}^-(IP(\Lambda^{op})) \times \mathcal{D}^-(IP(\Lambda)) \\
&\xrightarrow{\hat{\otimes}^L_\Lambda} \mathcal{D}^-(IP(R)) \\
&\xrightarrow{LH^{-n}} \mathcal{LH}(IP(R)),
\end{align*}
where in both cases the unlabelled maps are the obvious inclusions of full subcategories. Because $LH^n$ and $RH^n$ are cohomological functors, we get the usual long exact sequences in $\mathcal{LH}(IP(R))$ and $\mathcal{RH}(PD(R))$ coming from strict short exact sequences (in the appropriate category) in either variable, natural in both variables -- since these give distinguished triangles in the corresponding derived category.

When $\Ext_\Lambda^n$ or $\Tor^\Lambda_n$ take coefficients in $IP(\Lambda)$ or $PD(\Lambda)$, these coefficient modules should be thought of as objects in the appropriate left or right heart, via inclusion of the full subcategory.

\begin{rem}
\label{exttor}
The reason we cannot define `classical' derived functors in the sense of, say, \cite{Weibel} just in terms of topological module categories is essentially that these categories, like most interesting categories of topological modules, fail to be abelian. Intuitively this means that the naive definition of the homology of a chain complex of such modules -- that is, defining $$H_n(M) = \coker(\coim(d_{n-1}) \rightarrow \ker(d_n))$$ -- loses too much information. There is no well-behaved homology functor from chain complexes in a quasi-abelian category back to the category itself, so that a naive approach here fails. That is why we must use the more sophisticated machinery of passing to the left or right hearts, which function as `completions' of the original category to an abelian category, in an appropriate sense: see \cite{Schneiders}.

Our $\Ext$ and $\Tor$ functors are the appropriate analogues, in this setting, of the classical derived functors, in a sense made precise in the following proposition.
\end{rem}

Recall from Proposition \ref{de=dlhe} that we have equivalences $\mathcal{D}^-(IP(\Lambda)) \rightarrow \mathcal{D}^-(\mathcal{LH}(IP(\Lambda)))$ and $\mathcal{D}^+(PD(\Lambda)) \rightarrow \mathcal{D}^+(\mathcal{LH}(PD(\Lambda)))$. So we may think of $R\Hom_\Lambda^T(-,-)$ as a functor $$\mathcal{D}^-(\mathcal{LH}(IP(\Lambda))) \times \mathcal{D}^+(\mathcal{LH}(PD(\Lambda))) \rightarrow \mathcal{D}^+(\mathcal{LH}(PD(R)))$$ via these equivalences, and similarly for $\hat{\otimes}^L_\Lambda$. For the next proposition, we use these definitions.

Note that, by Remark \ref{ext0hom} below, $\Ext_\Lambda^0(-,-) \neq \Hom_\Lambda^T(-,-)$ as functors on $IP(\Lambda) \times PD(\Lambda)$.

\begin{prop}
$R\Hom_\Lambda^T(-,-)$ and $\Ext_\Lambda^n(-,-)$ are, respectively, the total derived functor and the $n$th classical derived functor of $\Ext_\Lambda^0(-,-)$. Similarly $\hat{\otimes}^L_\Lambda$ and $\Tor^\Lambda_n(-,-)$ are, respectively, the total derived functor and the $n$th classical derived functor of $\Tor^\Lambda_0(-,-)$.
\end{prop}
\begin{proof}
We prove the first statement; the second can be shown similarly. Write $$R\Ext_\Lambda^0(-,-): \mathcal{D}^-(\mathcal{LH}(IP(\Lambda))) \times \mathcal{D}^+(\mathcal{LH}(PD(\Lambda))) \rightarrow \mathcal{D}^+(\mathcal{LH}(PD(R)))$$ for the total right derived functor of $\Ext_\Lambda^0(-,-)$. Then for $M \in \mathcal{D}^-(\mathcal{LH}(IP(\Lambda)))$ with projective resolution $P$ and $N \in \mathcal{D}^+(\mathcal{LH}(PD(\Lambda)))$ with injective resolution $I$, $R\Ext_\Lambda^0(M,N)$ is by definition the total complex of the bicomplex $(\Ext_\Lambda^0(P_p,I^q))_{p,q}$. But $\Ext_\Lambda^0(P_p,I^q) = \Hom_\Lambda^T(P_p,I^q)$ because $P_p$ is projective and so is a resolution of itself. So the bicomplex is $(\Hom_\Lambda^T(P_p,I^q))_{p,q}$, and its total complex by definition is $R\Hom_\Lambda^T(M,N)$, giving the result for total derived functors. Taking $M \in \mathcal{LH}(IP(\Lambda))$ and $N \in \mathcal{LH}(PD(\Lambda))$, we get that the $n$th classical derived functor is $RH^n \circ R\Ext_\Lambda^0(M,N) = RH^n \circ R\Hom_\Lambda^T(M,N) = \Ext_\Lambda^n(M,N)$.
\end{proof}

\begin{lem}
\phantomsection
\label{rhom*}
\begin{enumerate}[(i)]
\item $R\Hom_\Lambda^T(-,-)$ and $\hat{\otimes}^L_\Lambda$ are Pontryagin dual in the sense that, given $M \in \mathcal{D}^-(IP(\Lambda))$ and $N \in \mathcal{D}^+(PD(\Lambda))$, there holds $R\Hom_\Lambda^T(M,N)^\ast = N^\ast \hat{\otimes}^L_\Lambda M$, naturally in 
$M,N$.
\item For $M \in \mathcal{LH}(IP(\Lambda))$ and $N \in \mathcal{RH}(PD(\Lambda))$, $\Ext_\Lambda^n(M,N)^\ast = \Tor^\Lambda_n(N^\ast,M)$.
\end{enumerate}
\end{lem}
\begin{proof}
We prove (i); (ii) follows by taking cohomology. Take a projective resolution $P$ of $M$ and an injective resolution $I$ of $N$, so that by duality $I^\ast$ is a projective resolution of $N^\ast$. Then $$R\Hom_\Lambda^T(M,N)^\ast = \Hom_\Lambda^T(P,I)^\ast = I^\ast \hat{\otimes}_\Lambda P = N^\ast \hat{\otimes}^L_\Lambda M,$$ naturally by the universal property of derived functors.
\end{proof}

\begin{rem}
More generally, as functors on the appropriate categories of bimodules, it follows from Theorem \ref{adjoint} that $L \hat{\otimes}^L_\Lambda -$ is left adjoint to the functor $R\Hom_\Theta^T(L,-)$ for $L \in IP(\Theta-\Lambda)$, and similarly for the $\Ext$ and $\Tor$ functors. Details are left to the reader.
\end{rem}

\begin{prop}
\phantomsection
\label{NastMast2}
\begin{enumerate}[(i)]
\item $R\Hom_\Lambda^T(M,N) = R\Hom_{\Lambda^{op}}^T(M^\ast,N^\ast)$ and 
$\Ext_\Lambda^n(M,N) = \Ext_{\Lambda^{op}}^n(N^\ast,M^\ast)$;
\item  $N^\ast \hat{\otimes}^L_\Lambda M = M \hat{\otimes}^L_{\Lambda^{op}} N^\ast$ and 
$\Tor^\Lambda_n(N^\ast,M) = \Tor^{\Lambda^{op}}_n(M,N^\ast)$;
\end{enumerate}
naturally in $M,N$.
\end{prop}
\begin{proof}
(ii) follows from (i) by Pontryagin duality. To see (i), take a projective resolution $P$ of $M$ and an injective resolution $I$ of $N$. Then $$R\Hom_\Lambda^T(M,N) = \Hom_\Lambda^T(P,I) = \Hom_{\Lambda^{op}}^T(I^\ast,P^\ast) = R\Hom_{\Lambda^{op}}^T(M^\ast,N^\ast),$$ by Lemma \ref{NastMast}. The rest follows by applying $LH^{-n}$.
\end{proof}

\begin{prop}
\label{1var}
$R\Hom_\Lambda^T$, $\Ext$, $\hat{\otimes}^L_\Lambda$ and $\Tor$ can be calculated using a resolution of either variable. That is, given $M$ with a projective resolution $P$ and $N$ with an injective resolution $I$, in the appropriate categories,
\begin{align*}
&R\Hom_\Lambda^T(M,N) = \Hom_\Lambda^T(P,N) = \Hom_\Lambda^T(M,I), \\
&\Ext_\Lambda^n(M,N) = RH^n(\Hom_\Lambda^T(P,N)) = H^n(\Hom_\Lambda^T(M,I)), \\
&N^\ast \hat{\otimes}^L_\Lambda M = N^\ast \hat{\otimes}_\Lambda P = I^\ast \hat{\otimes}_\Lambda M \text{ and} \\
&\Tor^\Lambda_n(N^\ast,M) = LH^{-n}(N^\ast \hat{\otimes}_\Lambda P) = LH^{-n}(I^\ast \hat{\otimes}_\Lambda M).
\end{align*}
\end{prop}
\begin{proof}
By Proposition \ref{bifunctor}, $R\Hom_\Lambda^T(M,N) = \Hom_\Lambda^T(M,I)$; everything else follows by some combination of Proposition \ref{NastMast2}, taking cohomology and applying Pontryagin duality.
\end{proof}

We will now see that, for module-theoretic purposes, it is sometimes more useful to apply $LH^n$ to right derived functors and $RH^n$ to left derived functors; though, as noted in Remark \ref{rrnotlr}, the resulting cohomological functors are not so well-behaved.

\begin{lem}
\label{hom=ext}
$LH^0 \circ R\Hom_\Lambda^T(M,N) = \Hom_\Lambda^T(M,N)$ and $RH^0 (N^\ast \hat{\otimes}^L_\Lambda M) = N^\ast \hat{\otimes}_\Lambda M$ for all $M \in IP(\Lambda), N \in PD(\Lambda)$, naturally in $M,N$.
\end{lem}
\begin{proof}
Take a projective resolution $P$ of $M$. Then
\begin{align*}
LH^0 \circ R\Hom_\Lambda^T(M,N) &= \ker(\Hom_\Lambda^T(P_0,N) \rightarrow \Hom_\Lambda^T(P_1,N)) \\
&= \Hom_\Lambda^T(\ker(P_0 \rightarrow P_1),N) \\
&= \Hom_\Lambda^T(M,N),
\end{align*}
because $\Hom_\Lambda^T$ commutes with kernels. The rest follows by duality.
\end{proof}

\begin{example}
$\mathbb{Z}_p$ is projective in $IP(\hat{\mathbb{Z}})$ by Corollary \ref{proj=tf}. Now consider the sequence $$0 \rightarrow \bigoplus_\mathbb{N} \mathbb{Z}_p \xrightarrow{f} \bigoplus_\mathbb{N} \mathbb{Z}_p \xrightarrow{g} \mathbb{Q}_p \rightarrow 0,$$ where $f$ is given by $(x_0,x_1,x_2,\ldots) \mapsto (x_0,x_1-p\cdot x_0,x_2-p\cdot x_1,\ldots)$ and $g$ is given by $(x_0,x_1,x_2,\ldots) \mapsto x_0 + x_1/p + x_2/p^2 + \cdots$. This sequence is exact on the underlying modules, so by Proposition \ref{ipsection} it is strict exact, and hence it is a projective resolution of $\mathbb{Q}_p$. By applying Pontryagin duality, we also get an injective resolution $$0 \rightarrow \mathbb{Q}_p \rightarrow \prod_\mathbb{N} \mathbb{Q}_p/\mathbb{Z}_p \rightarrow \prod_\mathbb{N} \mathbb{Q}_p/\mathbb{Z}_p \rightarrow 0.$$
\end{example}

Recall that, by Remark \ref{Qpnotinjective}, $\mathbb{Q}_p$ is not projective or injective.

\begin{lem}
\label{Qp0}
For all $n>0$ and all $M \in IP(\hat{\mathbb{Z}})$,
\begin{enumerate}[(i)]
\item $\Ext_{\hat{\mathbb{Z}}}^n(\mathbb{Q}_p,M^\ast) = 0$;
\item $\Ext_{\hat{\mathbb{Z}}}^n(M,\mathbb{Q}_p) = 0$;
\item $\Tor^{\hat{\mathbb{Z}}}_n(\mathbb{Q}_p,M) = 0$;
\item $\Tor^{\hat{\mathbb{Z}}}_n(M^\ast,\mathbb{Q}_p) = 0$.
\end{enumerate}
\end{lem}
\begin{proof}
By Lemma \ref{rhom*} and Proposition \ref{NastMast2}, it is enough to prove (iii). Since $\mathbb{Q}_p$ has a projective resolution of length $1$, the statement is clear for $n>1$. Now $\Tor^{\hat{\mathbb{Z}}}_1(\mathbb{Q}_p,M) = \ker(f \hat{\otimes}_{\hat{\mathbb{Z}}} M)$, in the notation of the example. Writing $M_p$ for $\mathbb{Z}_p \hat{\otimes}_{\hat{\mathbb{Z}}} M$, $f \hat{\otimes}_{\hat{\mathbb{Z}}} M$ is given by $$\bigoplus_\mathbb{N} M_p \rightarrow \bigoplus_\mathbb{N} M_p, (x_0,x_1,x_2,\ldots) \mapsto (x_0,x_1-p\cdot x_0,x_2-p\cdot x_1,\ldots),$$ because $\hat{\otimes}_{\hat{\mathbb{Z}}}$ commutes with direct sums. But this map is clearly injective, as required.
\end{proof}

\begin{rem}
\label{ext0hom}
By Lemma \ref{Qp0}, $\Ext_{\hat{\mathbb{Z}}}^0(\mathbb{Q}_p,-)$ is an exact functor from the category $\mathcal{RH}(PD(\hat{\mathbb{Z}}))$ to itself. In particular, writing $\mathcal{I}$ for the inclusion functor $PD(\hat{\mathbb{Z}}) \rightarrow \mathcal{RH}(PD(\hat{\mathbb{Z}}))$, the composite $\Ext_{\hat{\mathbb{Z}}}^0(\mathbb{Q}_p,-) \circ \mathcal{I}$ sends short strict exact sequences in $PD(\hat{\mathbb{Z}})$ to short exact sequences in $\mathcal{RH}(PD(\hat{\mathbb{Z}}))$ by Proposition \ref{reflective}. On the other hand, by Proposition \ref{reflective} again, the composite $\mathcal{I} \circ \Hom_{\hat{\mathbb{Z}}}^T(\mathbb{Q}_p,-)$ does not send short strict exact sequences in $PD(\hat{\mathbb{Z}})$ to short exact sequences in $\mathcal{RH}(PD(\hat{\mathbb{Z}}))$. Therefore, by \cite[Proposition 1.3.10]{Schneiders}, and in the terminology of \cite{Schneiders}, $\Hom_{\hat{\mathbb{Z}}}^T(\mathbb{Q}_p,-)$ is not \emph{RR left exact}: there is some short strict exact sequence $$0 \rightarrow L \rightarrow M \rightarrow N \rightarrow 0$$ in $PD(\hat{\mathbb{Z}})$ such that the induced map $\Hom_{\hat{\mathbb{Z}}}^T(\mathbb{Q}_p,M) \rightarrow \Hom_{\hat{\mathbb{Z}}}^T(\mathbb{Q}_p,N)$ is not strict. By duality, a similar result holds for tensor products with $\mathbb{Q}_p$.
\end{rem}

\section{Homology and cohomology of profinite\\ groups}

Let $G$ be a profinite group. We define the category of ind-profinite right $G$-modules, $IP(G^{op})$, to have as its objects ind-profinite abelian groups $M$ with a continuous map $M \times_k G \rightarrow M$, and as its morphisms continuous group homomorphisms which are compatible with the $G$-action. We define the category of pro-discrete $G$-modules, $PD(G)$, to have as its objects prodiscrete $\hat{\mathbb{Z}}$-modules $M$ with a continuous map $G \times M \rightarrow M$, and as its morphisms continuous group homomorphisms which are compatible with the $G$-action.

From now on, $R$ will denote a commutative profinite ring.

\begin{prop}
\phantomsection
\label{R+G=RG}
\begin{enumerate}[(i)]
\item $IP(G^{op})$ and $IP(\hat{\mathbb{Z}}\llbracket G \rrbracket^{op})$ are equivalent.
\item An ind-profinite right $R\llbracket G \rrbracket$-module is the same as an ind-profinite $R$-module $M$ with a continuous map $M \times_k G \rightarrow M$ such that $(mr)g = (mg)r$ for all $g \in G, r \in R, m \in M$.
\item $PD(G)$ and $PD(\hat{\mathbb{Z}}\llbracket G \rrbracket)$ are equivalent.
\item A pro-discrete $R\llbracket G \rrbracket$-module is the same as a pro-discrete $R$-module $M$ with a continuous map $G \times M \rightarrow M$ such that $g(rm) = r(gm)$ for all $g \in G, r \in R, m \in M$.
\end{enumerate}
\end{prop}
\begin{proof}
\begin{enumerate}[(i)]
\item Given $M \in IP(G^{op})$, take a cofinal sequence $\{M_i\}$ for $M$ as an ind-profinite abelian group. Replacing each $M_i$ with $M'_i = M_i \cdot G$ if necessary, we have a cofinal sequence for $M$ consisting of profinite right $G$-modules. By \cite[Proposition 5.3.6(c)]{R-Z}, each $M'_i$ canonically has the structure of a profinite right $\hat{\mathbb{Z}}\llbracket G \rrbracket$-module, and with this structure the cofinal sequence $\{M'_i\}$ makes $M$ into an object in $IP(\hat{\mathbb{Z}}\llbracket G \rrbracket^{op})$. This gives a functor $IP(G^{op}) \rightarrow IP(\hat{\mathbb{Z}}\llbracket G \rrbracket^{op})$. Similarly, we get a functor $IP(\hat{\mathbb{Z}}\llbracket G \rrbracket^{op}) \rightarrow IP(G^{op})$ by taking cofinal sequences and forgetting the $\hat{\mathbb{Z}}$-structure on the profinite elements in the sequence. These functors are clearly inverse to each other.
\item Similarly.
\item Similarly, replacing \cite[Proposition 5.3.6(c)]{R-Z} with \cite[Proposition 5.3.6(e)]{R-Z}.
\item Similarly.
\end{enumerate}
\end{proof}

By (ii) of Proposition \ref{R+G=RG}, given $M \in IP(R)$, we can think of $M$ as an object in $IP(R\llbracket G \rrbracket^{op})$ with trivial $G$-action. This gives a functor, the \emph{trivial module functor}, $IP(R) \rightarrow IP(R\llbracket G \rrbracket^{op})$, which clearly preserves strict exact sequences.

Given $M \in IP(R\llbracket G \rrbracket^{op})$, define the \emph{coinvariant module} $M_G$ by $$M/\overline{\langle m \cdot g - m, \text{ for all } g \in G, m \in M\rangle}.$$ This makes $M_G$ into an object in $IP(R)$. In the same way as for abstract modules, $M_G$ is the maximal quotient module of $M$ with trivial $G$-action, and so $-_G$ becomes a functor $IP(R\llbracket G \rrbracket^{op}) \rightarrow IP(R)$ which is left adjoint to the trivial module functor. We can define $-_G$ similarly for left ind-profinite $R\llbracket G \rrbracket$-modules.

By (iv) of Proposition \ref{R+G=RG}, given $M \in PD(R)$, we can think of $M$ as an object in $PD(R\llbracket G \rrbracket)$ with trivial $G$-action. This gives a functor which we also call the trivial module functor, $PD(R) \rightarrow PD(R\llbracket G \rrbracket)$, which clearly preserves strict exact sequences.

Given $M \in PD(R\llbracket G \rrbracket)$, define the \emph{invariant submodule} $M^G$ by $$\{m \in M: g\cdot m = m, \text{ for all } g \in G, m \in M\}.$$ It is a closed submodule of $M$, because $$M^G = \bigcap_{g \in G} \ker(M \rightarrow M, m \mapsto g\cdot m - m).$$ Therefore we can think of $M^G$ as an object in $PD(R)$. In the same way as for abstract modules, $M^G$ is the maximal submodule of $M$ with trivial $G$-action, and so $-^G$ becomes a functor $PD(R\llbracket G \rrbracket) \rightarrow PD(R)$ which is right adjoint to the trivial module functor. We can define $-^G$ similarly for right pro-discrete $R\llbracket G \rrbracket$-modules.

\begin{lem}
\phantomsection
\label{_G^G}
\begin{enumerate}[(i)]
\item For $M \in IP(R\llbracket G \rrbracket^{op})$, $M_G = M \hat{\otimes}_{R\llbracket G \rrbracket} R$.
\item For $M \in PD(R\llbracket G \rrbracket)$, $M^G = \Hom_{R\llbracket G \rrbracket}^T(R,M)$.
\end{enumerate}
\end{lem}
\begin{proof}
\begin{enumerate}[(i)]
\item Pick a cofinal sequence $\{M_i\}$ for $M$. By \cite[Lemma 6.3.3]{R-Z}, $(M_i)_G = M_i \hat{\otimes}_{R\llbracket G \rrbracket} R$, naturally in $M_i$. As a left adjoint, $-_G$ commutes with direct limits, so $$M_G = \varinjlim (M_i)_G = \varinjlim (M_i \hat{\otimes}_{R\llbracket G \rrbracket} R) = M \hat{\otimes}_{R\llbracket G \rrbracket} R$$ by Proposition \ref{tensorprops}.
\item Similarly, by \cite[Lemma 6.2.1]{R-Z}, because $-^G$ and $\Hom_{R\llbracket G \rrbracket}^T(R,-)$ commute with inverse limits.
\end{enumerate}
\end{proof}

\begin{cor}
Given $M \in IP(R\llbracket G \rrbracket^{op})$, $(M_G)^\ast = (M^\ast)^G$.
\end{cor}
\begin{proof}
Lemma \ref{_G^G} and Corollary \ref{adjointast}.
\end{proof}

We now define the \emph{$n$th homology} functor of $G$ over $R$ by $$H^R_n(G,-) = \Tor^{R\llbracket G \rrbracket}_n(-,R): \mathcal{LH}(IP(R\llbracket G \rrbracket^{op})) \rightarrow \mathcal{LH}(IP(R))$$ and the \emph{$n$th cohomology} functor of $G$ over $R$ by $$H_R^n(G,-) = \Ext_{R\llbracket G \rrbracket}^n(R,-): \mathcal{RH}(PD(R\llbracket G \rrbracket)) \rightarrow \mathcal{RH}(PD(R)).$$ As noted in Remark \ref{exttor}, we can also think of $H^R_n(G,-)$ as a functor $$IP(R\llbracket G \rrbracket^{op}) \rightarrow \mathcal{LH}(IP(R))$$ by precomposing with inclusion from these subcategories, and we may do so without further comment.

We have by Lemma \ref{rhom*} that:

\begin{prop}
$H^R_n(G,M)^\ast = H_R^n(G,M^\ast)$ for all $M \in \mathcal{LH}(IP(R\llbracket G \rrbracket^{op}))$, naturally in $M$.
\end{prop}

Of course, one can calculate all these objects using the projective resolution of $R$ arising from the usual bar resolution, \cite[Section 6.2]{R-Z}, and this shows that the homology and cohomology are unchanged if we forget the $R$-module structure and think of $M$ as an object of $\mathcal{LH}(IP(\hat{\mathbb{Z}}\llbracket G \rrbracket^{op}))$; that is, the underlying complex of abelian $k$-groups of $H^R_n(G,M)$, and the underlying complex of topological abelian groups of $H_R^n(G,M^\ast)$, are $H^{\hat{\mathbb{Z}}}_n(G,M)$ and $H_{\hat{\mathbb{Z}}}^n(G,M^\ast)$, respectively. We therefore write
\begin{align*}
H_n(G,M) &= H^{\hat{\mathbb{Z}}}_n(G,M) \text{ and} \\
H^n(G,M^\ast) &= H_{\hat{\mathbb{Z}}}^n(G,M^\ast).
\end{align*}

\begin{thm}[Universal Coefficient Theorem]
Suppose $M \in PD(\hat{\mathbb{Z}}\llbracket G \rrbracket)$ has trivial $G$-action. Then there are non-canonically split short strict exact sequences $$0 \rightarrow \Ext_{\hat{\mathbb{Z}}}^1(H_{n-1}(G,\hat{\mathbb{Z}}),M) \rightarrow H^n(G,M) \rightarrow \Ext_{\hat{\mathbb{Z}}}^0(H_n(G,\hat{\mathbb{Z}}),M) \rightarrow 0,$$ $$0 \rightarrow \Tor^{\hat{\mathbb{Z}}}_0(M^\ast,H_n(G,\hat{\mathbb{Z}})) \rightarrow H_n(G,M^\ast) \rightarrow \Tor^{\hat{\mathbb{Z}}}_1(M^\ast,H_{n-1}(G,\hat{\mathbb{Z}})) \rightarrow 0.$$
\end{thm}
\begin{proof}
We prove the first sequence; the second follows by Pontryagin duality. Take a projective resolution $P$ of $\hat{\mathbb{Z}}$ in $IP(\hat{\mathbb{Z}}\llbracket G \rrbracket)$ with each $P_n$ profinite, so that $H^n(G,M) = H^n(\Hom_{\hat{\mathbb{Z}}\llbracket G \rrbracket}^T(P,M))$. Because $M$ has trivial $G$-action, $M = \Hom_{\hat{\mathbb{Z}}}^T(\hat{\mathbb{Z}},M)$, where we think of $\hat{\mathbb{Z}}$ as an ind-profinite $\hat{\mathbb{Z}}-\hat{\mathbb{Z}}\llbracket G \rrbracket$-bimodule with trivial $G$-action. So
\begin{align*}
\Hom_{\hat{\mathbb{Z}}\llbracket G \rrbracket}^T(P,M) &= \Hom_{\hat{\mathbb{Z}}\llbracket G \rrbracket}^T(P,\Hom_{\hat{\mathbb{Z}}}^T(\hat{\mathbb{Z}},M)) \\
&= \Hom_{\hat{\mathbb{Z}}}^T(\hat{\mathbb{Z}} \hat{\otimes}_{\hat{\mathbb{Z}}\llbracket G \rrbracket} P,M) \\
&= \Hom_{\hat{\mathbb{Z}}}^T(P_G,M).
\end{align*}
Note that $P_G$ is a complex of profinite modules, so all the maps involved are automatically strict. Since $-_G$ is left adjoint to an exact functor (the trivial module functor), we get in the same way as for abelian categories that $-_G$ preserves projectives, so each $(P_n)_G$ is projective in $IP(\hat{\mathbb{Z}})$ and hence torsion free by Corollary \ref{proj=tf}. Now the profinite subgroups of each $(P_n)_G$ consisting of cycles and boundaries are torsion-free and hence projective in $IP(\hat{\mathbb{Z}})$ by Corollary \ref{proj=tf}, so $P_G$ splits. Then the result follows by the same proof as in the abstract case, \cite[Section 3.6]{Weibel}.
\end{proof}

\begin{cor}
For all $n$,
\begin{enumerate}[(i)]
\item $H_n(G,\mathbb{Z}_p) = \Tor^{\hat{\mathbb{Z}}}_0(\mathbb{Z}_p,H_n(G,\hat{\mathbb{Z}})) = \mathbb{Z}_p \hat{\otimes}_{\hat{\mathbb{Z}}} H_n(G,\hat{\mathbb{Z}})$.
\item $H_n(G,\mathbb{Q}_p) = \Tor^{\hat{\mathbb{Z}}}_0(\mathbb{Q}_p,H_n(G,\hat{\mathbb{Z}}))$.
\item $H^n(G,\mathbb{Q}_p) = \Ext_{\hat{\mathbb{Z}}}^0(H_n(G,\hat{\mathbb{Z}}),\mathbb{Q}_p)$.
\end{enumerate}
\end{cor}
\begin{proof}
(i) holds because $\mathbb{Z}_p$ is projective; (ii) and (iii) follow from Lemma \ref{Qp0}.
\end{proof}

Suppose now that $H$ is a (profinite) subgroup of $G$. We can think of $R \llbracket G \rrbracket$ as an ind-profinite $R\llbracket H \rrbracket - R\llbracket G \rrbracket$-bimodule: the left $H$-action is given by left multiplication by $H$ on $G$, and the right $G$-action is given by right multiplication by $G$ on $G$. We will denote this bimodule by $R \llbracket ^{H\searrow}G^{\swarrow G} \rrbracket$.

If $M \in IP(R\llbracket G \rrbracket)$, we can restrict the $G$-action to an $H$-action. Moreover, maps of $G$-modules which are compatible with the $G$-action are compatible with the $H$-action. So restriction gives a functor $$\Res^G_H: IP(R\llbracket G \rrbracket) \rightarrow IP(R\llbracket H \rrbracket).$$ $\Res^G_H$ can equivalently be defined by the functor $R \llbracket ^{H\searrow}G^{\swarrow G} \rrbracket \hat{\otimes}_{R \llbracket G \rrbracket} -$. Similarly, we can define a restriction functor $$\Res^G_H: PD(R\llbracket G \rrbracket^{op}) \rightarrow PD(R\llbracket H \rrbracket^{op})$$ by $\Hom_{R \llbracket G \rrbracket}^T(R \llbracket ^{H\searrow}G^{\swarrow G} \rrbracket,-)$.

On the other hand, given $M \in IP(R\llbracket H \rrbracket^{op})$, $M \hat{\otimes}_{R\llbracket H \rrbracket} R \llbracket ^{H\searrow}G^{\swarrow G} \rrbracket$ becomes an object in $IP(R\llbracket G \rrbracket^{op})$. In this way, $- \hat{\otimes}_{R\llbracket H \rrbracket} R \llbracket ^{H\searrow}G^{\swarrow G} \rrbracket$ becomes a functor, induction, $$\Ind^G_H: IP(R\llbracket H \rrbracket^{op}) \rightarrow IP(R\llbracket G \rrbracket^{op}).$$ Also, $\Hom_{R \llbracket H \rrbracket}^T(R \llbracket ^{H\searrow}G^{\swarrow G} \rrbracket,-)$ becomes a functor, coinduction, which we denote by $$\Coind^G_H: PD(R\llbracket H \rrbracket) \rightarrow PD(R\llbracket G \rrbracket).$$

Since $R \llbracket ^{H\searrow}G^{\swarrow G} \rrbracket$ is projective in $IP(R\llbracket H \rrbracket)$ and $IP(R\llbracket G \rrbracket)^{op}$, $\Res^G_H, \Ind^G_H$ and $\Coind^G_H$ all preserve strict exact sequences. Moreover, $\Res^G_H$ and $\Ind^G_H$ commute with colimits of ind-profinite modules because tensor products do, and $\Res^G_H$ and $\Coind^G_H$ commute with limits of pro-discrete modules because $\Hom$ does in the second variable.

We can similarly define restriction on right ind-profinite or left pro-discrete $R\llbracket G \rrbracket$-modules, induction on left ind-profinite $R\llbracket G \rrbracket$-modules and coinduction on right pro-discrete $R\llbracket G \rrbracket$-modules, using $R \llbracket ^{G\searrow}G^{\swarrow H} \rrbracket$. Details are left to the reader.

Suppose an abelian group $M$ has a left $H$-action together with a topology that makes it into both an ind-profinite $H$-module and a pro-discrete $H$-module. For example, this is the case if $M$ is second-countable profinite or countable discrete. Then both $\Ind^G_H$ and $\Coind^G_H$ are defined. When $H$ is open in $G$, we get $\Ind^G_H- = \Coind^G_H-$ in the same way as the abstract case, \cite[Lemma 6.3.4]{Weibel}.

\begin{lem}
\label{indast}
For $M \in IP(R\llbracket H \rrbracket^{op})$, $(\Ind^G_HM)^\ast = \Coind^G_H(M^\ast)$. For $N \in IP(R\llbracket G \rrbracket^{op})$, $(\Res^G_HN)^\ast = \Res^G_H(N^\ast)$.
\end{lem}
\begin{proof}
\begin{align*}
(\Ind^G_HM)^\ast &= (M \hat{\otimes}_{R\llbracket H \rrbracket} R \llbracket ^{H\searrow}G^{\swarrow G} \rrbracket)^\ast \\
&= \Hom_{R \llbracket H \rrbracket}^T(R \llbracket ^{H\searrow}G^{\swarrow G} \rrbracket,M^\ast) = \Coind^G_H(M^\ast).
\end{align*}
\begin{align*}
(\Res^G_HN)^\ast &= (N \hat{\otimes}_{R \llbracket G \rrbracket} R \llbracket ^{G\searrow}G^{\swarrow H} \rrbracket)^\ast \\
&= \Hom_{R \llbracket H \rrbracket}^T(R \llbracket ^{G\searrow}G^{\swarrow H} \rrbracket, N^\ast) = \Res^G_H(N^\ast).
\end{align*}
\end{proof}

\begin{lem}
\phantomsection
\label{indadjoint}
\begin{enumerate}[(i)]
\item $\Ind^G_H$ is left adjoint to $\Res^G_H$. That is, for $M \in IP(R\llbracket H \rrbracket)$, $N \in IP(R\llbracket G \rrbracket)$, $\Hom_{R\llbracket G \rrbracket}^{IP}(\Ind^G_HM,N) = \Hom_{R\llbracket H \rrbracket}^{IP}(M,\Res^G_HN)$, naturally in $M$ and $N$.
\item $\Coind^G_H$ is right adjoint to $\Res^G_H$. That is, for $M \in PD(R\llbracket G \rrbracket)$, $N \in PD(R\llbracket H \rrbracket)$, $\Hom_{R\llbracket G \rrbracket}^{PD}(M,\Coind^G_HN) = \Hom_{R\llbracket H \rrbracket}^{PD}(\Res^G_HM,N)$, naturally in $M$ and $N$.
\end{enumerate}
\end{lem}
\begin{proof}
(i) and (ii) are equivalent by Pontryagin duality and Lemma \ref{indast}. We show (i). Pick cofinal sequences $\{M_i\},\{N_j\}$ for $M,N$. Then
\begin{align*}
\Hom_{R\llbracket G \rrbracket}^{IP}(\Ind^G_HM,N) &= \Hom_{R\llbracket G \rrbracket}^{IP}(\varinjlim (\Ind^G_H M_i),\varinjlim N_j) \\
&= \varprojlim_i\varinjlim_j\Hom_{R\llbracket G \rrbracket}^{IP}(\Ind^G_H M_i,N_j) \\
&= \varprojlim_i\varinjlim_j\Hom_{R\llbracket H \rrbracket}^{IP}(M_i,\Res^G_HN_j) \\
&= \Hom_{R\llbracket H \rrbracket}^{IP}(\varinjlim M_i,\varinjlim \Res^G_HN_j) \\
&= \Hom_{R\llbracket H \rrbracket}^{IP}(M,\Res^G_HN)
\end{align*}
by Lemma \ref{ipprojinj} and the Pontryagin dual of \cite[Lemma 6.10.2]{R-Z}, and all the isomorphisms in this sequence are natural.
\end{proof}

\begin{cor}
\label{indproj}The functor
$\Ind^G_H$ sends projectives in $IP(R\llbracket H \rrbracket)$ to projectives in $IP(R\llbracket G \rrbracket)$. Dually, $\Coind^G_H$ sends injectives in $PD(R\llbracket H \rrbracket)$ to injectives in $PD(R\llbracket G \rrbracket)$.
\end{cor}
\begin{proof}
The adjunction of Lemma \ref{indadjoint} shows that, for $P \in IP(R\llbracket H \rrbracket)$ projective, $\Hom_{R\llbracket G \rrbracket}^{IP}(\Ind^G_HP,-) = \Hom_{R\llbracket H \rrbracket}^{IP}(P,\Res^G_H-)$ sends strict epimorphisms to surjections, as required.

The second statement follows from the first by applying the result for $\Ind^G_H$ to $IP(R\llbracket H \rrbracket^{op})$, and then using Pontryagin duality.
\end{proof}

\begin{lem}
\label{indres}
For $M \in IP(R\llbracket H \rrbracket^{op})$, $N \in IP(R\llbracket G \rrbracket)$, $\Ind^G_HM \hat{\otimes}_{R\llbracket G \rrbracket} N = M \hat{\otimes}_{R\llbracket H \rrbracket} \Res^G_H N$ and $\Hom_{R \llbracket G \rrbracket}^T(N,\Coind^G_H(M^\ast)) = \Hom_{R \llbracket H \rrbracket}^T(N,M^\ast)$, naturally in $M,N$.
\end{lem}
\begin{proof}
$\Ind^G_HM \hat{\otimes}_{R\llbracket G \rrbracket} N = M \hat{\otimes}_{R\llbracket H \rrbracket} R \llbracket ^{H\searrow}G^{\swarrow G} \rrbracket \hat{\otimes}_{R\llbracket G \rrbracket} N = M \hat{\otimes}_{R\llbracket H \rrbracket} \Res^G_H N$. The second equation follows by applying Pontryagin duality and Lemma \ref{indast}.
\end{proof}

\begin{thm}[Shapiro's Lemma]
For $M \in \mathcal{D}^-(IP(R\llbracket H \rrbracket^{op}))$, $N \in \mathcal{D}^-(IP(R\llbracket G \rrbracket))$, we have:
\begin{enumerate}[(i)]
\item $\Ind^G_HM \hat{\otimes}_{R\llbracket G \rrbracket}^L N = M \hat{\otimes}_{R\llbracket H \rrbracket}^L \Res^G_HN$;
\item $N \hat{\otimes}_{R\llbracket G \rrbracket^{op}}^L \Ind^G_HM = \Res^G_HN \hat{\otimes}_{R\llbracket H \rrbracket^{op}}^L M$;
\item $R\Hom_{R\llbracket G \rrbracket^{op}}^T(\Ind^G_HM,N^\ast) = R\Hom_{R\llbracket H \rrbracket^{op}}^T(M,\Res^G_HN^\ast)$;
\item $R\Hom_{R\llbracket G \rrbracket}^T(N,\Coind^G_HM^\ast) = R\Hom_{R\llbracket H \rrbracket}^T(\Res^G_HN,M^\ast)$;
\end{enumerate}
naturally in $M,N$. Similar statements hold for the $\Ext$ and $\Tor$ functors.
\end{thm}
\begin{proof}
We show (i); (ii)-(iv) follow by Lemma \ref{rhom*} and Proposition \ref{NastMast2}. Take a projective resolution $P$ of $M$. By Corollary \ref{indproj}, $\Ind^G_HP$ is a projective resolution of $\Ind^G_HM$. Then
\begin{align*}
\Ind^G_HM \hat{\otimes}_{R\llbracket G \rrbracket}^L N &= \Ind^G_HP \hat{\otimes}_{R\llbracket G \rrbracket} N \\
&= P \hat{\otimes}_{R\llbracket H \rrbracket} \Res^G_HN \text{ by Lemma \ref{indres}} \\
&= M \hat{\otimes}_{R\llbracket H \rrbracket}^L \Res^G_HN,
\end{align*}
and all these isomorphisms are natural. For the rest, apply the cohomology functors.
\end{proof}

\begin{cor}
For $M \in IP(R\llbracket H \rrbracket^{op})$,
\begin{align*}
H^R_n(G,\Ind^G_HM) &= H^R_n(H,M) \text{ and} \\
H_R^n(G,\Coind^G_HM^\ast) &= H_R^n(H,M^\ast)
\end{align*}
for all $n$, naturally in $M$.
\end{cor}
\begin{proof}
Apply Shapiro's Lemma with $N=R$ with trivial $G$-action -- the restriction of this action to $H$ is also trivial.
\end{proof}

If $K$ is a profinite normal subgroup of $G$, then for $M \in IP(R \llbracket G \rrbracket^{op})$, $M_K$ becomes an
ind-profinite right $R \llbracket G/K \rrbracket$-module, as in the abstract case. So we may think of $-_K$ as a functor $IP(R \llbracket G \rrbracket^{op}) \rightarrow IP(R \llbracket G/K \rrbracket^{op})$ and consider its right derived functor $$R(-_K): \mathcal{D}^-(IP(R\llbracket G \rrbracket^{op})) \rightarrow \mathcal{D}^-(IP(R\llbracket G/K \rrbracket^{op}));$$ we write $H^R_s(K,-)$ for the `classical' derived functor given by the composition
\begin{align*}
\mathcal{LH}(IP(R\llbracket G \rrbracket^{op})) &\rightarrow \mathcal{D}^-(IP(R\llbracket G \rrbracket^{op})) \\
&\xrightarrow{R(-_K)} \mathcal{D}^-(IP(R\llbracket G/K \rrbracket^{op})) \\
&\xrightarrow{\mathcal{LH}^{-s}} \mathcal{LH}(IP(R\llbracket G/K \rrbracket^{op})).
\end{align*}
Thus we can compose the two functors $H^R_s(K,-)$ and $H^R_r(G/K,-)$.

The case of $-^K$ can be handled similarly.

\begin{thm}[Lyndon-Hochschild-Serre Spectral Sequence]
Suppose $K$ is a profinite normal subgroup of $G$. Then there are bounded spectral sequences $$E^2_{rs} = H^R_r(G/K,H^R_s(K,M)) \Rightarrow H^R_{r+s}(G,M)$$ for all $M \in \mathcal{LH}(IP(R\llbracket G \rrbracket^{op}))$ and $$E_2^{rs} = H_R^r(G/K,H_R^s(K,M)) \Rightarrow H_R^{r+s}(G,M)$$ for all $M \in \mathcal{RH}(PD(R\llbracket G \rrbracket))$, both naturally in $M$. In particular, these hold for $M \in IP(R\llbracket G \rrbracket^{op})$ and $M \in PD(R \llbracket G \rrbracket)$, respectively.
\end{thm}
\begin{proof}
We prove the first statement; then Pontryagin duality gives the second by Lemma \ref{rhom*}. By the universal properties of $-_K$, $-_{G/K}$ and $-_G$, it is easy to see that $(-_K)_{G/K} = -_G$. Moreover, as for abstract modules, $-_K$ is left adjoint to the forgetful functor $IP(R\llbracket G/K \rrbracket^{op}) \rightarrow IP(R\llbracket G \rrbracket^{op})$, which sends strict exact sequences to strict exact sequences, and hence $-_K$ preserves projectives. So the result is just an application of the Grothendieck Spectral Sequence, Theorem \ref{grothendieck}.
\end{proof}

\begin{rem}
One must be careful in applying this spectral sequence: no such spectral sequence exists in general if we try to define derived functors back to the original module categories, for the reasons discussed in Remark \ref{exttor}. The naive definition of homology functor is not sufficiently well-behaved here.
\end{rem}

\section{Comparison to other cohomologies}
\label{compare}

Let $P(\Lambda)$ and $D(\Lambda)$ be the categories of profinite and discrete $\Lambda$-modules, both with continuous homomorphisms. We will think of $P(\Lambda)$ as a full subcategory of $IP(\Lambda)$ and $D(\Lambda)$ as a full subcategory of $PD(\Lambda)$. We consider alternative definitions of $\Ext_\Lambda^n$ using these categories, and show how they compare to our definition. Specifically, we will compare our definitions to:
\begin{enumerate}[(i)]
\item the classical cohomology of profinite rings using discrete coefficients, found for instance in \cite{R-Z};
\item the theory of cohomology for profinite modules of type $\FP_\infty$ over profinite rings, developed in \cite{S-W}, allowing profinite coefficients;
\item the continuous cochain cohomology, defined as in \cite{Tate}, for all topological modules over topological rings;
\item the reduced continuous cochain cohomology, defined as in \cite{CG}, for all topological modules over topological rings.
\end{enumerate}

Recall from Section \ref{qacs} that the inclusion $\mathcal{I}^{op}: PD(\Lambda) \rightarrow \mathcal{RH}(PD(\Lambda))$ has a right adjoint $\mathcal{C}^{op}$. We can give an explicit description of these functors by duality: for $M \in PD(\Lambda)$, $\mathcal{I}^{op}(M) = (0 \rightarrow M \rightarrow 0)$. Each object in $\mathcal{RH}(PD(\Lambda))$ is isomorphic to a complex $M' = (0 \rightarrow M^0 \xrightarrow{f} M^1 \rightarrow 0)$ in $PD(\Lambda)$, where $M^0$ is in degree $0$ and $f$ is epic, and $\mathcal{C}^{op}(M') = \ker(f)$. Also the functors $$RH^n: \mathcal{D}(PD(\Lambda)) \rightarrow \mathcal{RH}(PD(\Lambda))$$ are given by $$RH^n(\cdots \xrightarrow{d^{n-1}} M^n \xrightarrow{d^n} M^{n+1} \xrightarrow{d^{n+1}} \cdots) = (0 \rightarrow \coker(d^{n-1}) \rightarrow \im(d^n) \rightarrow 0),$$ with $\coker(d^{n-1})$ in degree $0$.

Given $M \in P(\Lambda), N \in D(\Lambda)$, to avoid ambiguity we write $\underline{\Hom}_\Lambda(M,N)$ for the discrete $R$-module of continuous $\Lambda$-homomorphisms $M \rightarrow N$; we have $\underline{\Hom}_\Lambda(M,N) = \Hom^T_\Lambda(M,N)$ in this case. Let $P$ be a projective resolution of $M$ in $P(\Lambda)$ and $I$ an injective resolution of $N$ in $D(\Lambda)$: recall that projectives in $P(\Lambda)$ are projective in $IP(\Lambda)$ and injectives in $D(\Lambda)$ are injective in $PD(\Lambda)$ by Lemma \ref{pdinjective}. In \cite{R-Z}, the derived functors of $$\underline{\Hom}_\Lambda: P(\Lambda) \times D(\Lambda) \rightarrow D(R)$$ are defined by $$\underline{\Ext}_\Lambda^n(M,N) = H^n(\underline{\Hom}_\Lambda(P,N)),$$ or equivalently by $H^n(\underline{\Hom}_\Lambda(M,I))$, where cohomology is taken in $D(R)$.

\begin{prop}
\label{underlineExt}
$\mathcal{I}^{op}\underline{\Ext}_\Lambda^n(M,N) = \Ext_\Lambda^n(M,N)$ as pro-discrete $R$-modules.
\end{prop}
\begin{proof}
We have $\Ext_\Lambda^n(M,N) = RH^n(\Hom^T_\Lambda(P,N))$. Because each $P_n$ is profinite, $\Hom^T_\Lambda(P,N) = \underline{\Hom}_\Lambda(P,N)$ is a cochain complex of discrete $R$-modules; write $d^n$ for the map $\Hom^T_\Lambda(P_n,N) \rightarrow \Hom^T_\Lambda(P_{n+1},N)$. In the abelian category $D(R)$, applying the Snake Lemma to the diagram
\[
\xymatrix{& \im(d^{n-1}) \ar[r] \ar[d] & \Hom^T_\Lambda(P_n,N) \ar[r] \ar[d] & \coker(d^{n-1}) \ar[r] \ar[d] & 0 \\
0 \ar[r] & \ker(d^n) \ar[r] & \Hom^T_\Lambda(P_n,N) \ar[r] & \coim(d^n)}
\]
shows that
\begin{align*}
H^n(\underline{\Hom}_\Lambda(P,N)) &= \coker(\im(d^{n-1}) \rightarrow \ker(d^n)) \\
&= \ker(\coker(d^{n-1}) \rightarrow \coim(d^n)).
\end{align*}
Next, using once again that $D(R)$ is abelian, we have
\begin{align*}
RH^n(\Hom^T_\Lambda(P,N)) &= (0 \rightarrow \coker(d^{n-1}) \rightarrow \im(d^n) \rightarrow 0) \\
&= (0 \rightarrow \coker(d^{n-1}) \rightarrow \coim(d^n) \rightarrow 0),
\end{align*}
so it is enough to show that the map of complexes
\[
\xymatrix{0 \ar[r] \ar[d] & \ker(\coker(d^{n-1}) \rightarrow \coim(d^n)) \ar[r] \ar[d] & 0 \ar[r] \ar[d] & 0 \ar[d] \\
0 \ar[r] & \coker(d^{n-1}) \ar[r] & \coim(d^n) \ar[r] & 0}
\]
is a strict quasi-isomorphism, or equivalently that its cone $$0 \rightarrow \ker(\coker(d^{n-1}) \rightarrow \coim(d^n)) \rightarrow \coker(d^{n-1}) \rightarrow \coim(d^n) \rightarrow 0$$ is strict exact, which is clear.
\end{proof}

If on the other hand we are given $M,N \in P(\Lambda)$ with $M$ finitely generated, we avoid ambiguity by writing $\Hom^P_\Lambda(M,N)$ for the profinite $R$-module (with the compact-open topology) of continuous $\Lambda$-homomorphisms $M \rightarrow N$. Then, writing $P(\Lambda)_\infty$ for the full subcategory of $P(\Lambda)$ whose objects are of type $\FP_\infty$, in \cite{S-W} the derived functors of $$\Hom^P_\Lambda: P(\Lambda)_\infty \times P(\Lambda) \rightarrow P(R)$$ are defined by $$\Ext_\Lambda^{P,n}(M,N) = H^n(\Hom^P_\Lambda(P,N)),$$ where $P$ is a projective resolution of $M$ in $P(\Lambda)$ such that each $P_n$ is finitely generated, and cohomology is taken in $P(R)$. Assume that $N$ is second-countable, so that $\Ext_\Lambda^n(M,N)$ is defined. Because $P(R)$ is an abelian category, the same proof as Proposition \ref{underlineExt} shows:

\begin{prop}
$\mathcal{I}^{op}\Ext_\Lambda^{P,n}(M,N) = \Ext_\Lambda^n(M,N)$ as pro-discrete $R$-modules.
\end{prop}

For any $M \in IP(\Lambda)$ and $N \in T(\Lambda)$, the $R$-module of continuous $\Lambda$-homomorphisms $\operatorname{cHom}_\Lambda(M,N)$, with the compact-open topology, defines a functor $IP(\Lambda) \times T(\Lambda) \rightarrow TAb$, where $TAb$ is the category of topological abelian groups and continuous homomorphisms. For a projective resolution $P$ of $M$ in $IP(\Lambda)$, the \emph{continuous cochain $\Ext$ functors} are then defined by $$\cExt_\Lambda^n(M,N) = H^n(\operatorname{cHom}_\Lambda(P,N)),$$ 
where the cohomology is taken in $TAb$. That is, $$\cExt_\Lambda^n(M,N) = \ker(d^n)/\coim(d^{n-1}),$$ where $\ker(d^n)$ is given the subspace topology and $\ker(d^n)/\coim(d^{n-1})$ is given the quotient topology. For $\Lambda=\hat{\mathbb Z}\llbracket G \rrbracket$ and $M=\hat{\mathbb Z}$ the trivial $\Lambda$-module, this definition essentially coincides with the continuous cochain cohomology of $G$ introduced in \cite{Tate}[Section 2], and the results stated here for $\Ext$ functors easily translate to the special case of group cohomology, which we leave to the reader. Indeed, it is easy to check that the bar resolution described in \cite{Tate} gives a projective resolution of $\hat{\mathbb{Z}}$ in $IP(\hat{\mathbb Z}\llbracket G \rrbracket)$, and hence that the cohomology theory described there coincides with ours.

Given a short exact sequence $$0 \rightarrow A \rightarrow B \rightarrow C \rightarrow 0$$ of topological $\Lambda$-modules, we do not in general get a long exact sequence of $\cExt$ functors. If the short exact sequence is such that the sequence of underlying modules of $$0 \rightarrow \operatorname{cHom}_\Lambda(P_n,A) \rightarrow \operatorname{cHom}_\Lambda(P_n,B) \rightarrow \operatorname{cHom}_\Lambda(P_n,C) \rightarrow 0$$ is exact for all $P_n$, then (by forgetting the topology) we do get a sequence $$0 \rightarrow \cExt_\Lambda^0(M,A) \rightarrow \cExt_\Lambda^0(M,B) \rightarrow \cExt_\Lambda^0(M,C) \rightarrow \cdots$$ which is a long exact sequence of the underlying modules.

In general, we cannot expect $\cExt_\Lambda^n(M,N)$ to be a Hausdorff topological group, since the images of the continuous homomorphisms $$d^{n-1}: \operatorname{cHom}_\Lambda(P_{n-1},N)\rightarrow \operatorname{cHom}_\Lambda(P_n,N)$$ are not necessarily closed. This immediately suggests the following alternative definition. We define the \emph{reduced continuous cochain $\Ext$ functors} by $$\rExt_\Lambda^n(M,N) = \ker(d^n)/\overline{\coim(d^{n-1})},$$ with the quotient topology. Clearly $\rExt$ coincides with $\cExt$ exactly when the set-theoretic image of $\coim(d^{n-1})$ is closed in $\ker(d^n)$. Note that, even when we have a long exact sequence in $\cExt$, the passage to $\rExt$ need not be exact.

For $M \in IP(\Lambda)$ and $N \in PD(\Lambda)$, let $P$ be a projective resolution of $M$ in $IP(\Lambda)$ 
and $(\Hom_\Lambda^T(P,N),d)$ be the associated cochain complex. Then 
$$\Ext_\Lambda^n(M,N) = (0 \rightarrow \coker(d^{n-1}) \xrightarrow{f} \im(d^n) \rightarrow 0).$$

\begin{prop}
In this notation, 
\begin{enumerate}[(i)]
\item $\rExt_\Lambda^n(M,N) = \ker(f)$. 
\item If $d^{n-1}(\Hom_\Lambda^T(P_{n-1},N))$ is closed in 
$\Hom_\Lambda^T(P_n,N)$, then there holds $\cExt_\Lambda^n(M,N) = \ker(f)$.
\end{enumerate}
\end{prop}
\begin{proof}
Consider the diagram
\[
\xymatrix{0 \ar[r] & \im(d^{n-1}) \ar[r] \ar[d] & \Hom_\Lambda^T(P_n,N) \ar[r] \ar[d]^{=} & \coker(d^{n-1}) \ar[r] \ar[d] & 0 \\
0 \ar[r] & \ker(d^n) \ar[r] & \Hom_\Lambda^T(P_n,N) \ar[r] & \coim(d^n) \ar[r] & 0,}
\]
with the obvious maps. The rows are strict exact, and the vertical maps are clearly strict, so after applying Pontryagin duality, Lemma \ref{snake} says that $$\ker(\coker(d^{n-1}) \rightarrow \coim(d^n)) \cong \coker(\im(d^{n-1}) \rightarrow \ker(d^n)) = \rExt_\Lambda^n(M,N).$$ On the other hand, $$\ker(\coker(d^{n-1}) \rightarrow \coim(d^n)) = \ker(\coker(d^{n-1}) \rightarrow \im(d^n))$$ because $\coim(d^n) \rightarrow \im(d^n)$ is monic. The second statement is clear.
\end{proof}

\begin{cor}
\label{cext}
$\cExt_\Lambda^n(M,N) = \ker f$ for all $n$ in the following two cases:
\begin{enumerate}[(i)]
\item $M \in P(\Lambda), N \in D(\Lambda)$;
\item $M,N \in P(\Lambda)$ with $M$ of type $\FP_\infty$ and $N$ second-countable.
\end{enumerate}
\end{cor}
\begin{proof}
In these cases, $\Hom_\Lambda^T(P,N)$ is in the abelian categories $D(\Lambda)$ and $P(\Lambda)$ respectively, so it is strict and the conditions for the proposition are satisfied, for all $n$.
\end{proof}

On the other hand, for any $M \in IP(\Lambda)$ and $N \in PD(\Lambda)$, we can also consider the alternative cohomological functors mentioned in Remark \ref{rrnotlr}: $$LH^n \circ R\Hom_\Lambda^T(M,N) = (0 \rightarrow \coim(d^{n-1}) \xrightarrow{g} \ker(d^n) \rightarrow 0).$$ Then we can recover the continuous cochain $\Ext$ functors from this information:

\begin{prop}
\begin{enumerate}[(i)]
\item $\cExt_\Lambda^n(M,N) = \coker(g)$, where the cokernel is taken in $T(R)$;
\item $\rExt_\Lambda^n(M,N) = \coker(g)$, where the cokernel is taken in $PD(R)$.
\end{enumerate}
\end{prop}
\begin{proof}
\begin{enumerate}[(i)]
\item In $T(R)$, there holds $\coker(g) = \ker(d^n)/\coim(d^{n-1})$, with the quotient topology, which is $\cExt_\Lambda^n(M,N)$ by definition.
\item Similarly, in $PD(R)$, $\coker(g) = \ker(d^n)/\overline{\coim(d^{n-1})}$.
\end{enumerate}
\end{proof}

From another perspective, this proposition says that all the $\Ext$ functors we have considered can be obtained from the total derived functor $R\Hom_\Lambda^T(-,-)$.

Exactly the same approach as this section makes it possible to compare our $\Tor$ functors to other definitions, with similar conclusions. We leave both definitions and proofs to the reader, noting only the following results. For $M \in IP(\Lambda)$ and $N \in IP(\Lambda^{op})$, let $P$ be a projective resolution of $M$ in $IP(\Lambda)$ 
and $(N \hat{\otimes}_\Lambda P,d)$ be the associated chain complex. Then 
$$\Tor^\Lambda_n(N,M) = (0 \rightarrow \coim(d^{n-1}) \xrightarrow{h} \ker(d^n) \rightarrow 0).$$

\begin{prop}
\begin{enumerate}[(i)]
\item The \emph{continuous chain $\Tor$ functor} $\cTor^\Lambda_n(M,N)$ (defined in the obvious way) is the cokernel $\coker(h)$, taken in $T(R)$;
\item the \emph{reduced continuous chain $\Tor$ functor} $\rTor^\Lambda_n(M,N)$ (defined in the obvious way) is the cokernel $\coker(h)$, taken in $IP(R)$;
\item $\cTor^\Lambda_n(M,N) = \rTor^\Lambda_n(M,N)$ if and only if $h(\coim(d^{n-1}))$ is closed in $\ker(d^n)$.
\end{enumerate}
\end{prop}

By Pontryagin duality and Corollary \ref{cext}, the condition for (iii) is satisfied for all $n$ if $M \in P(\Lambda)$, $N \in P(\Lambda^{op})$, or if $M \in P(\Lambda)$ is of type $\FP_\infty$ and $N \in D(\Lambda^{op})$ is discrete and countable.

\end{document}